\documentclass[preprint,aap,11pt]{persocls}

\startlocaldefs
\endlocaldefs
\usepackage{amsmath,amssymb,latexsym,psfrag}
\usepackage{amsfonts}
\usepackage{graphicx}
\usepackage{epsf,epsfig,color}

\newcommand{\sone}{s_2}
\newcommand{\stwo}{s_1}

\newtheorem{theo}{Theorem}

\newtheorem{prop}{Proposition}
\newtheorem{lemm}{Lemma}
\newtheorem{defn}{Definition}

\def\E{\mathbb{E}}
\def\S{\mathbb{S}}
\def\0{{\bf 0}}
\def\Z{\mathbb{Z}}
\def\R{\mathbb{R}}
\def\PP{\mathbb{P}}

\def\B{{ B}}

\def\en2{{\epsilon_n^2}}
\renewcommand{\E}{\mathbb E \,}
\newcommand{\C}{{\cal C}}

\newcommand{\tod}{\stackrel{{\cal D}}{\longrightarrow}}
\newcommand{\eqd}{\stackrel{{\cal D}}{=}}

\newcommand{\eqco}{\setcounter{equation}{0}}
\newcommand{\thco}{\setcounter{theo}{0}}
\newcommand{\prco}{\setcounter{prop}{0}}
\newcommand{\laco}{\setcounter{lemm}{0}}
\newcommand{\coco}{\setcounter{coro}{0}}
\newcommand{\cjco}{\setcounter{conj}{0}}

\newcommand{\deco}{\setcounter{defn}{0}}
\newcommand{\allco}{\eqco  \thco \prco \laco \coco \cjco \deco}
\newcommand{\X}{{\cal X}}

\def\2e{{\epsilon^2_\la}}
\def\la{{\lambda}}

\def\v{{ w }}

\newcommand{\un}{{ h^{(\la, \alpha)} }}
\newcommand{\ul}{{ h^{(\lambda, \alpha)} }}

\newcommand{\Y}{{\cal Y}}

\renewcommand{\P}{{{\cal P}}}

\newcommand{\Var}{{\rm Var}}

\newcommand{\Vol}{{\rm Vol}}

\newcommand{\F}{{\cal F}}

\def\bdm{\begin{displaymath}}
\newcommand{\edm}{\end{displaymath}}
\def\benu{\begin{enumerate}}
\def\eenu{\end{enumerate}}
\def\beqn{\begin{equation}}
\def\eeqn{\end{equation}}
\def\be{\begin{equation}}
\def\ee{\end{equation}}
\def\bea{\begin{eqnarray}}
\def\eea{\end{eqnarray}}
\newcommand{\bean}{\begin{eqnarray*}}
\newcommand{\eean}{\end{eqnarray*}}
\newcommand{\bear}{\begin{eqnarray}}
\newcommand{\eear}{\end{eqnarray}}

\newcommand{\1}{{\bf 1}}

\def\R{\mathbb{R}}
\def\B^2{\mathbb{D}}
\def\S{\mathbb{S}}
\def\B{\mathbb{B}}

\def\qed{\hfill\hbox{${\vcenter{\vbox{
    \hrule height 0.4pt\hbox{\vrule width 0.4pt height 6pt
    \kern5pt\vrule width 0.4pt}\hrule height 0.4pt}}}$}}

\def\la{{\lambda}}

\begin{document}

\begin{frontmatter}

% "Title of the paper"
\title{Convex hulls of perturbed random point sets}
\runtitle{Convex hulls of perturbed random point sets}

% indicate corresponding author with \corref{}
% \author{\fnms{John} \snm{Smith}\corref{}\ead[label=e1]{smith@foo.com}\thanksref{t1}}
% \thankstext{t1}{Thanks to somebody}
% \address{line 1\\ line 2\\ printead{e1}}
% \affiliation{Some University}

\begin{aug}
\author{\fnms{Pierre} \snm{Calka}\thanksref{t1}\ead[label=e1]{pierre.calka@univ-rouen.fr}}
\and
\author{\fnms{J.~E.} \snm{Yukich}\thanksref{t2}\ead[label=e2]{jey0@lehigh.edu}}
\affiliation{Universit\'e de Rouen Normandie and Lehigh University}
\address{Pierre Calka\\Universit\'e de Rouen Normandie\\LMRS, avenue de l'universit\'e, BP 12\\F-76801 Saint-\'Etienne-du-Rouvray, FRANCE\\ \printead{e1}}
\address{J.~E. Yukich\\Lehigh University\\14 East Packer Ave, Christmas-Saucon\\Bethlehem, Pa. 18015, USA\\ \printead{e2}}
\thankstext{t1}{Research supported in part by the {\it Institut Universitaire de France}, the French ANR grant ASPAG (ANR-17-CE40-0017) and the Norman RIN grant ALENOR}
\thankstext{t2}{Research supported by NSF grant  DMS-1406410 and a Simons Collaboration Grant}
\runauthor{P. Calka and J.~E. Yukich}
\end{aug}

%\and
%\author{\fnms{???} \snm{???}\ead[label=e2]{???}}
%\address{\printead{e2}}
%\affiliation{???}

\begin{abstract}
We consider the convex hull of the perturbed point process comprised of $n$ i.i.d. points, each distributed as the sum of a uniform point on the unit sphere $\S^{d-1}$ and a uniform point in the $d$-dimensional ball centered at the origin and of radius $n^{\alpha}, \alpha \in (-\infty, \infty)$. This model, inspired by the smoothed complexity analysis introduced in computational geometry \cite{DGGT,ST}, is a perturbation of the classical random polytope. We show that the perturbed point process, after rescaling, converges in the scaling limit to one of five Poisson point processes according to whether $\alpha$ belongs to one of five regimes. The intensity measure of the
limit Poisson point process undergoes a transition at the values $\alpha = \frac{-2} {d -1}$ and $\alpha = \frac{2} {d + 1}$ and it gives rise to four  rescalings for the $k$-face functional on perturbed data.    These rescalings are used to establish explicit expectation asymptotics for the number of $k$-dimensional faces of the convex hull of either perturbed binomial or Poisson data.
In the case of Poisson input, we  establish explicit variance asymptotics and a central limit theorem for the number of $k$-dimensional faces.  Finally it is shown that the rescaled boundary of the convex hull of the perturbed point process converges to the boundary of a parabolic hull process.
\end{abstract}



\end{frontmatter}

\section{Introduction and main results}\label{INTRO}

\allco

%\noindent {\bf 1.1.  Introduction.}
The study of random polytopes has developed hand in hand with the construction of algorithms in computational geometry. The underlying idea is that the size of randomly generated geometric objects gives information on the complexity of an algorithm allowing the construction of these objects. In particular, the number of extreme points of the convex hull provides bounds on the running time of algorithms used to construct this hull.

On the computational geometry side, algorithms are generally analyzed in either the mean- or worst-case analyses. By mean-case analysis, we mean the expected time needed to construct the convex hull when the input is purely random while worst-case refers to the maximal time needed to make the construction when all possible inputs are considered.  In practice, neither of these two complexities well approximates the execution speed of a geometric algorithm. In particular, in their breakthrough paper \cite{ST}, Spielman and Teng chose the example of the simplex algorithm whose speed is often  polynomial in many application domains whereas the theory indicates that the complexity can be exponential. To explain this fact, they substitute a new complexity for the mean- and worst-case complexities, namely the smoothed complexity which interpolates between  them. It consists in perturbing an arbitrary point set  and calculating the maximum of the complexity over all possible initial configurations. The perturbation considered by Spielman and Teng is Gaussian. They claim that such smoothed complexity is more realistic since on the one hand the worst-case analysis is too pessimistic while  on the other hand the mean-case analysis is irrelevant because the considered data inherently possess a  structure making  it differ from randomly generated data.
Their main result is the polynomial growth of the smoothed complexity for the simplex algorithm; see Theorem 5.1 therein.

 A few years later, Devillers et al. \cite{DGGT} extended the notion of smoothed complexity to the convex hull construction. In this work, the smoothed complexity is the maximum over all possible initial configurations of the mean number of extreme points when the points from the initial configuration are randomly perturbed. Two possible perturbations are considered:  uniform distributions in Euclidean balls and isotropic  Gaussian distributions. Both models have one degree of freedom,  which is the radius of the perturbation ball for the former and the standard deviation of the Gaussian distribution for the latter. Sharp $\Theta$-estimates for the smoothed complexity in both cases are proved and in particular, several regimes are identified.

Inspired by the notion of smoothed complexity, we introduce a new random polytope model, one generated by perturbed input and which goes as follows. We do not calculate the   smoothed complexity itself but instead we take a particular random initial configuration comprised of $n$ points  independently and uniformly distributed  on the $d$-dimensional unit sphere $\S^{d-1}$. The points are all in convex position and thus we are initially in the worst-case set-up. Perturbations presumably bring one closer to the mean-case set-up.  We perturb each point by a random vector uniformly distributed in
the ball $B_d(\0, n^{\alpha})$, $\alpha \in \R$. Here $B_d(x,r)$ is the Euclidean ball centered at $x\in\R^d$ and of radius $r \in (0, \infty)$.
One expects that for a wide range of $\alpha$ such perturbed configurations should provide a number of extreme points whose average is close to the maximum over initial configurations.

 A second major motivation comes from convex geometry.  The classical study of random convex polytopes assumes that the polytopes are generated by data which is uncorrupted by noise, a somewhat unrealistic assumption. It is natural to ask about the effect of noise on the number of $k$-dimensional faces,  $k \in \{0,...,d-1 \}$,  in polytopes generated by an i.i.d. sample of size $n$. Up to now, there has been little research in this direction. When the noise components are uniformly distributed in $B_d(\0, n^{\alpha})$, $\alpha \in \R$, one may expect that small perturbations will not modify the number of $k$-dimensional faces, whereas large perturbations will have a substantial impact.  This paper addresses these and related questions.

 \vskip.5cm

 \noindent{\bf 1.1.  The polytope models.} We now describe the two models considered in this paper and our main contributions.

 \vskip.1cm

 \noindent{\em Binomial model}. Let $X_i, 1 \leq i \leq n$, be independently and  uniformly distributed on the $(d-1)$-dimensional unit sphere $\S^{d-1}, d \geq 2$.
Given $\alpha \in \R$, each $X_i, 1 \leq i \leq n$, is perturbed by a random vector $e_i:= e_{i,n}(\alpha)$, with $e_i, 1 \leq i \leq n,$ denoting a collection of i.i.d. random variables uniformly distributed on  $B_d(\0, n^{\alpha})$.  The $e_i, 1 \leq i \leq n,$ represent input errors or noise.
 Put $\tilde{X}_i:=  X_i + e_i$, $\X_n := \{X_i \}_{i=1}^n$, $\tilde{\X}_{n, \alpha} := \{ \tilde{X}_i \}_{i=1}^n$, and let $K_{n, \alpha}$ denote the convex hull of $\tilde{\X}_{n, \alpha}$.

\vskip.1cm

\noindent{\em Poisson model.}  For all $\la > 0$, let $\P_\la$ be a Poisson point process of intensity $\la/(d \kappa_d)$ on $\S^{d - 1}, d \geq 2$.
Here $\kappa_d$ denotes the volume of the unit ball in $\R^d$. Write $\P_\la:= \{ x_i \}_{i = 1}^{N(\la)}$, where $N(\la)$ is a Poisson random variable with parameter $\la$.  Given $\alpha \in \R$, each $x_i, 1 \leq i \leq N(\la)$, is perturbed by a random vector
$e_i:= e_{i, \la}(\alpha)$, with $e_i$ denoting an independent random variable uniformly distributed on $B_d(\0, \la^{\alpha})$.
%$e_i:= \{e_{i, \la}\}_{i = 1}^{N(\alpha)}$, a collection of i.i.d. random variables uniformly distributed on $B_d(\0, \la^{\alpha})$.
 Put $\tilde{x}_i:= \tilde{x}_{i, \lambda} (\alpha):= x_i + e_i$,  $\tilde{\P}_{\la, \alpha}:= \{ \tilde{x}_i \}_{i=1}^{N(\la)}$, and let $K_{\la,\alpha}$ denote the convex hull of $\tilde{\P}_{\la, \alpha}$.  Integral subscripts are reserved to denote binomial input, i.e., $K_{n, \alpha}$ denotes the convex hull of $\tilde{\X}_{n, \alpha}$.
 \vskip.1cm

\begin{figure} [!h]
\centering
\includegraphics[scale=0.6]{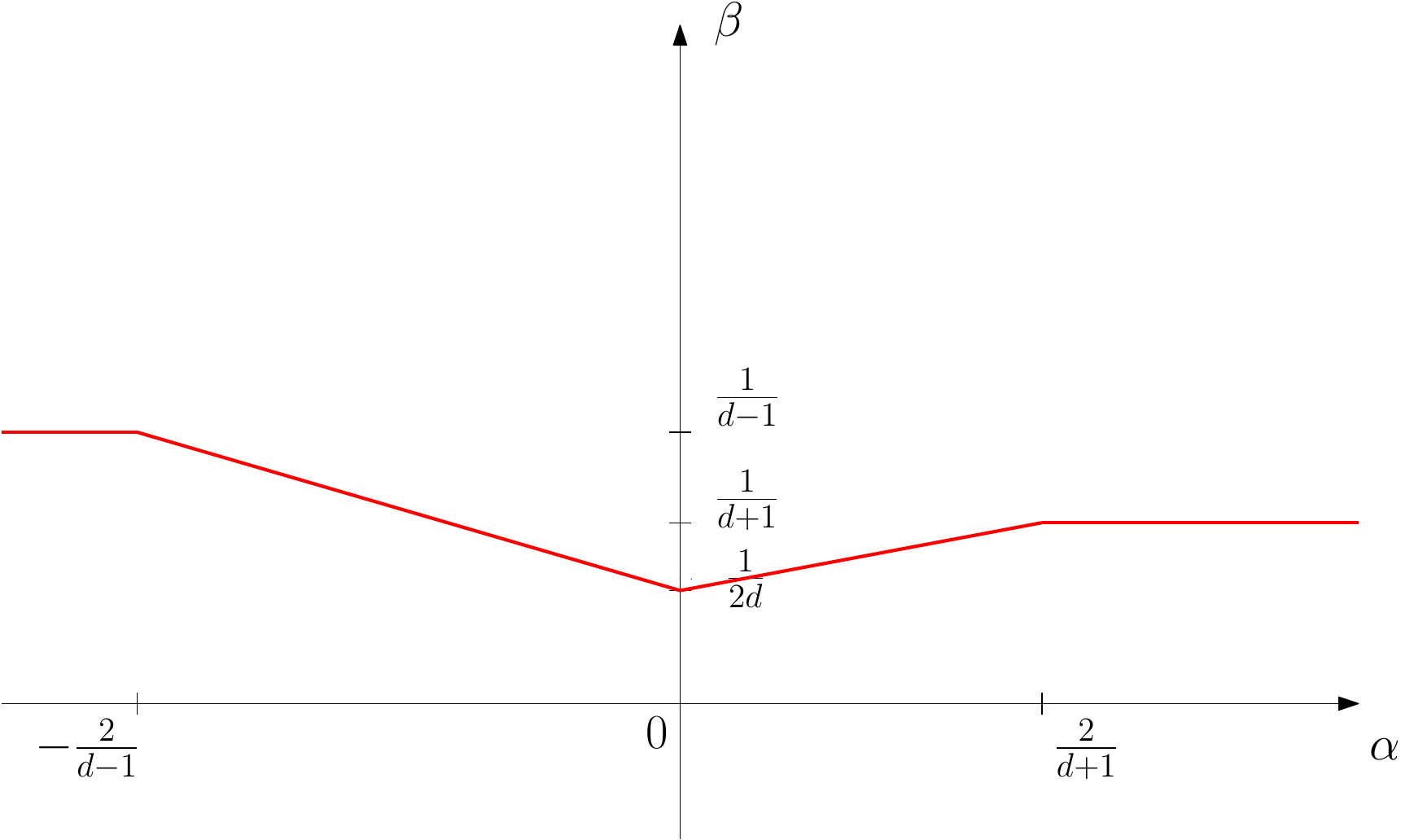}
\caption{The graph of $\beta(\alpha)$}
\label{fig:graphbeta}
\end{figure}

Our main contributions  are as follows:
\begin{itemize}

\item we show that, after rescaling, the perturbed point processes $\tilde{\X}_{n, \alpha}$ and  $\tilde{\P}_{\la, \alpha}$, converge as $n$ and $\la$ tend to  infinity, respectively,
to one of five scaling limit Poisson point processes $\P^{(\infty, \alpha) }$ {\em according to whether $\alpha$ belongs to one of these five regimes} : \ $( \frac{2}{d + 1} , \infty)$, $\{  \frac{2}{d + 1} \}$, $({-2 \over d -1}, \frac{2}{d + 1})$, $ \{ {-2 \over d -1} \} $, $(- \infty, {-2 \over d -1} )$.  The density of the intensity measure of the scaling limit $\P^{(\infty, \alpha)}$  undergoes a transition at  the  values $\alpha = \frac{2} {d + 1}$ and $\alpha = \frac{-2} {d -1}$.  At the value $\alpha = \frac{-2} {d -1}$, one observes a `phase transition' as the support of the measure changes from $\R^{d-1} \times \R^+$ to $\R^{d-1} \times \{\0\}.$

 %\item  There are four natural re-scalings for the $k$-face functional on perturbed data, which themselves undergo a sort of phase transition at the two %critical values and also at $\alpha = 0$.  \red{The expected $k$-face functional is not monotonically varying with $\alpha$, at least in the large $\la$ %regime.}

\item as in \cite{DGGT}, we show that the mean $k$-face functional displays four distinct rescalings on perturbed binomial and Poisson input, $k \in \{0,...,d-1 \}$. We go beyond \cite{DGGT} by providing  explicit limit constants for the  properly rescaled mean $k$-face functional for all  $\alpha$ including the critical $\alpha$ values equalling $\frac{2} {d + 1}$, $0$, and $ \frac{-2} {d -1}$, where there is a transition in the scaling regimes. The expected $k$-face functional is not monotonically varying with $\alpha$, at least in the large $\la$ regime; cf. Figure \ref{fig:graphbeta}.

\item we establish explicit variance asymptotics and rates of normal convergence for the rescaled $k$-face functional on perturbed Poisson input for all $\alpha$.

\item we show that the boundary of the convex hull of the perturbed point process converges, after rescaling, to the boundary of the parabolic hull process associated with $\P^{(\infty, \alpha) }$ on $\R^{d-1} \times \R^+$.

\end{itemize}

Our approach depends heavily on parabolic scaling transforms of the perturbed point set. When identifying the scaling limit point process for non-perturbed uniform data in the unit ball $B_d(\0,1)$, the authors previously showed \cite{CSY,CY, SY} that one may usefully introduce a parabolic  scaling transform.  A parabolic transform is also the key to determining the scaling limit of
non-perturbed data with an isotropic distribution \cite{GT1, GT2}.
Here we show that one needs  not one but six parabolic scaling transforms to properly identify the scaling limit Poisson point process (the regime  $(\frac{-2} {d - 1}, \frac{2} {d + 1} )$ requires three transforms, one for negative scalars, one for zero,  and one for positive scalars whereas the two regimes  $\{ \frac{-2} {d - 1} \}$ and $( \frac{2} {d +1}, \infty)$ require the same transform).  We are optimistic that the scaling transforms and techniques of this paper may be used to develop  the limit theory for i.i.d. samples  having an arbitrary isotropic distribution in $\R^d$ and which are  subject to perturbations more general than those considered here.  The case of Gaussian perturbations is the subject of a future paper.

\vskip.1cm

\noindent {\bf 1.2.  Scaling factors, scaling transforms and scaling limits.} Scaling exponents $\beta:= \beta(\alpha)$ are defined according to regimes for $\alpha$; see the first column of Table \ref{tab:summary} and the graph of $\beta(\cdot)$ in Figure \ref{fig:graphbeta}.

\begin{table}[!h]\label{tab:summary}
\begin{tabular}{|c|c|c|c|}
\hline
$\alpha$&$\beta:=\beta(\alpha)$&$u_{\la,\alpha}$&$\nu^{(\infty, \alpha)}$\\
\hline\hline
&&&\\
$\left(\frac{2}{d+1},\infty\right)$&$\frac{1}{d+1}$&$\kappa_d^{-\frac{1}{d+1}}\la^{\frac{1}{d+1}}$&$\mathrm{d}h$\\
&&&\\
\hline
&&&\\
$\frac{2}{d+1}$&$\frac{2+\alpha(d-1)}{4d}$&$\kappa_d^{-\frac{1}{d+1}}\la^{\frac{1}{d+1}}$&$\stwo(\kappa_d^{\frac{2}{d+1}}h)\mathrm{d}h$ \\
&&&\\
\hline
&&&\\
$\left(0,\frac{2}{d+1}\right)$&$\frac{2+\alpha(d-1)}{4d}$ &$ \left(\frac{2^{\frac{d-1}{2}}\kappa_{d-1}}{d\kappa_d^2}\right)^{\frac{1}{2d}}\la^{\beta}$&$h^{\frac{d-1}{2}}\mathrm{d}h$\\&&&\\
\hline
&&&\\$0$&$\frac{1}{2d}$&$\sqrt{2}\left(\frac{2^{\frac{d-1}{2}}\kappa_{d-1}}{d\kappa_d^2}\right)^{\frac{1}{2d}}\la^{\frac{1}{2d}}$&$h^{\frac{d-1}{2}}\mathrm{d}h$\\&&&\\\hline
&&&\\
$\left(-\frac{2}{d-1},0\right)$&$\frac{2-\alpha(d+1)}{4d}$&$ \left(\frac{2^{\frac{d-1}{2}}\kappa_{d-1}}{d\kappa_d^2}\right)^{\frac{1}{2d}}\la^{\beta}$&$h^{\frac{d-1}{2}}\mathrm{d}h$\\&&&\\\hline&&&\\
$-\frac{2}{d-1}$&$\frac{1}{d-1}$&$\kappa_d^{-\frac{1}{d+1}}\la^{\frac{1}{d-1}}$&$\sone(\kappa_d^{\frac{2}{d+1}}h)\mathrm{d}h$ \\&&&\\\hline&&&\\
$\left(-\infty,-\frac{2}{d-1}\right)$&$\frac{1}{d-1}$&$(d\kappa_d)^{-\frac{1}{d-1}}\la^{\frac{1}{d-1}}$&$\delta_{0}(h)$\\&&&\\
\hline
\end{tabular}
~\\~\\
%}
\caption{Summary of the different regimes}
\thispagestyle{empty}
\end{table}

These exponents give rise to scaling factors $u_{\la,\alpha}$
governing the growth of the $k$-face functional and which are central to all that follows; see the second column of Table \ref{tab:summary}.
There is thus a scaling factor assigned to each of the six (possibly degenerate) intervals:  $(- \infty, \frac{-2} {d-1} ), \{ \frac{-2} {d-1}  \} , ( \frac{-2} {d-1}, 0), \{ 0 \}, (0, \frac{2} {d + 1} )$ and $[ \frac{2} {d + 1}, \infty)$.

To each of these six regimes we associate a scaling
 transformation $T^{(\la, \alpha)}: \ B_d(\0, 1 + \la^{\alpha}) \to \R^{d-1} \times \R^+$ given by
\be \label{map2}
x \mapsto \left( u_{\la,\alpha} \exp_{d-1}^{-1} \frac{x} {|x|}, \ u_{\la,\alpha}^2 (1 - \frac{|x|}{1 + \la^{\alpha}} )\right).
\ee
Here $|x|$ denotes the Euclidean norm in $\R^d$ and the exponential map $\exp_{d-1}$ sends a point $v$ in the plane $T(\S^{d-1})$ which is tangent
to ${\bf u}_0:= (0,0,...,1) \in \R^d$ to a point on $\S^{d-1}$ lying at a geodesic distance $|v|$ from ${\bf u}_0$.  Here $|v|$ is the Euclidean norm in the tangent plane.
The radial scaling factor $\ u_{\la,\alpha}^2$ is the square of the angular scaling factor $u_{\la,\alpha}$.  When $\la \to \infty$ this parabolic scaling means that  $T^{(\la, \alpha)}$ carries balls and half-spaces to up- and down-paraboloids, respectively (cf. Section 3.1).

For all $h \geq 0$, let $\stwo(h)$ be the normalized surface area of the cap formed by intersecting $B_d({\bf u}_0,1)$ with $\R^{d-1} \times [0,h]$. Thus for $h \in [0,2]$ we have
\begin{align}  \label{eq:defs2}
\stwo(h) & :=\frac{1}{d\kappa_d}\sigma_{d-1}(\S^{d-1}\cap (\R^{d-1} \times [0,h])) \nonumber \\
& =\frac{(d-1)\kappa_{d-1}}{d\kappa_d}\int_0^{\arccos(1-h)}\sin^{d-2}\theta\mathrm{d}\theta,
\end{align}
where ${\sigma}_{d-1}$ is the $(d-1)$-dimensional Hausdorff measure of $\S^{d-1}$, and otherwise $\stwo(h) = 1$.

For all $h \geq 0$, let $\sone(h)$ be the surface area of the intersection of $B_d({\bf u}_0,1)$ with $\R^{d-1} \times \{h\}$ normalized
 by the surface area of $\S^{d-1}$.  Thus for  $h \in [0,2]$
\begin{equation}
  \label{eq:defs1}
\sone(h) :=\frac{\kappa_{d-1}}{d\kappa_d}(h(2-h))^{\frac{d-1}{2}}
\end{equation}
and otherwise $\sone(h) = 0$.
For $h$ small we have  $\sone(h) \sim \stwo(h) \sim \frac{\kappa_{d-1}}{d\kappa_d}(2h)^{(d-1)/2 }$.  Here and elsewhere we write $f(h) \sim g(h)$ to denote
$\lim f(h)/g(h) = 1$, whenever $h \to 0$ or $h \to \infty$, depending on context.

\begin{figure}[!h]
 \centering
\includegraphics[scale=0.4]{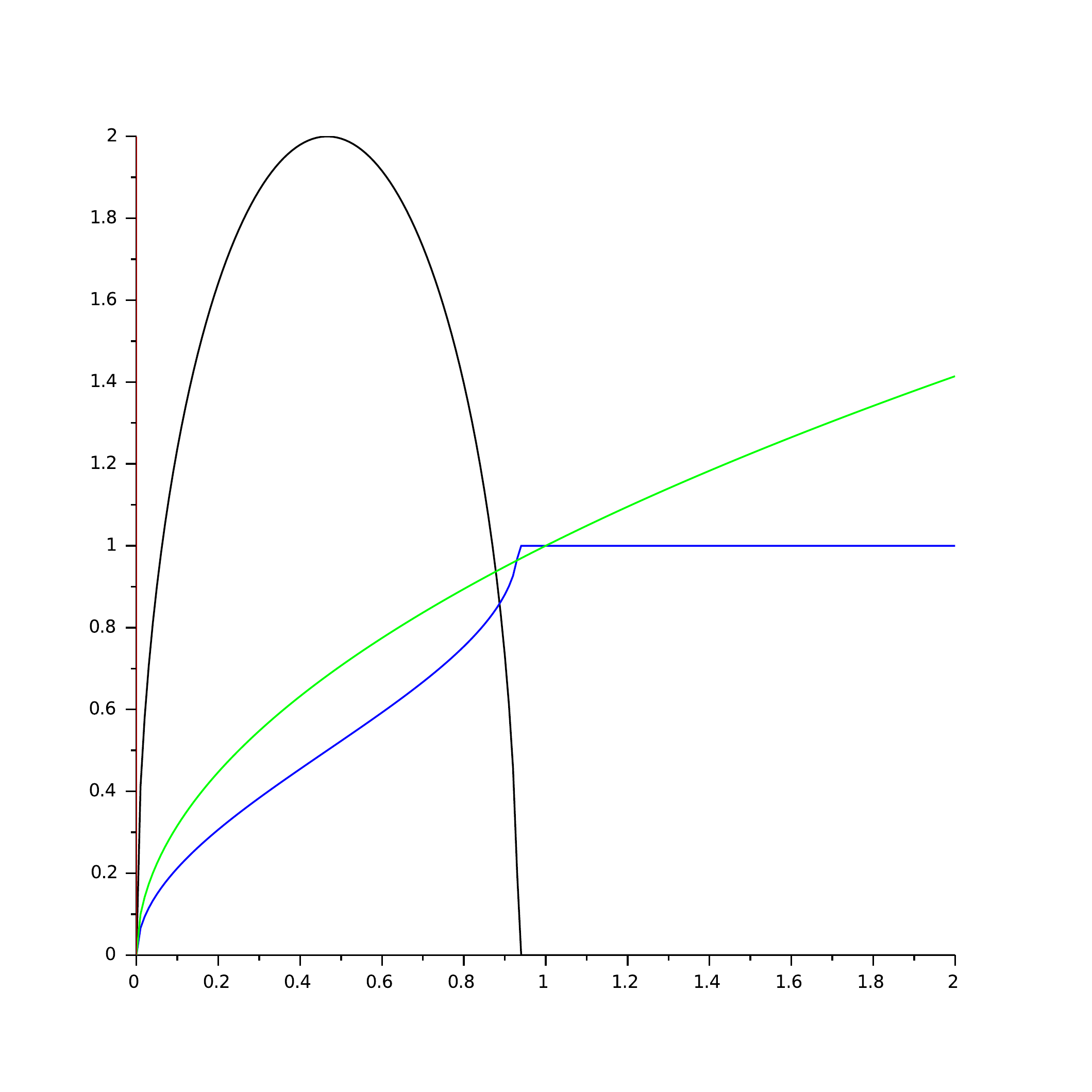}
 \caption{The measure $\nu^{(\infty,\alpha)}$ of ${\mathcal P}^{(\infty,\alpha)}$ for $d = 2$ and as a function of $h\ge 0$ for $\alpha=\frac{2}{d+1}$ (blue), $\alpha \in ( \frac{- 2} {d - 1}, \frac{2}{d+1})$ (green), $\alpha=-\frac{2}{d-1}$ (black)  and $\alpha\in (-\infty,-\frac{2}{d-1})$ (red)}
 \label{fig:densities}
\end{figure}

Let $\delta_x$ stand for the point measure putting mass one at $x$.
We let $\P^{(\infty, \alpha) }$ be the Poisson point process on $\R^{d-1} \times \R^+$ of intensity measure
\begin{equation}
  \label{eq:deflimitingintensity}
\mu^{(\infty, \alpha)}(\mathrm{d}v,\mathrm{d}h)=\mathrm{d}v \times \nu^{(\infty, \alpha)}(\mathrm{d}h)
\end{equation}
 where $\nu^{(\infty, \alpha)}(\mathrm{d}h)$ is given in the third column of Table \ref{tab:summary}; see also Figure \ref{fig:densities} for the graph of the density of $\nu^{(\infty,\alpha)}$.

The transformed perturbed binomial and Poisson point processes converge to one of five scaling limit point processes, according to the regimes for
 $\alpha$.

\begin{prop} \label{prop1} For $\alpha \in (-\infty, \infty)$ we have as $n \to \infty$
\be \label{Tnconv}
T^{(n, \alpha)} (\tilde{\X}_{n, \alpha}) \tod \P^{(\infty, \alpha) } .
\ee
We also have as $\la \to \infty$
\be \label{Tnconv1}
T^{(\la, \alpha)} (\tilde{\P}_{\la, \alpha}  ) \tod \P^{(\infty, \alpha) }.
\ee
\end{prop}

The point processes  $\P^{(\infty, \alpha)}, \alpha \in (-\infty, \infty)$, figure prominently in our main results, as seen in the next subsection.

\vskip.1cm

\noindent {\bf 1.3.  Limit theory for the $k$-face functional on perturbed input.} Let $\P$ be a Poisson point process on $\R^{d-1} \times \R$.
Let $\Pi^{\downarrow}:= \{(v,h) \in \R^{d-1} \times \R: \ h \leq \frac{|v|^2} {2} \}$ and put $\Pi^{\downarrow}(w):= w \oplus
\Pi^{\downarrow}$, where $\oplus$ denotes Minkowski addition.  Consider the maximal union, here denoted by
$\Phi(\P)$, of  parabolic grains  $\Pi^{\downarrow}(w), w \in \R^{d-1} \times \R$,
having the property that $\Pi^{\downarrow}(w)$ belongs to
$\Phi(\P)$ if its interior contains no point of $\P$.  Thus
%\be \label{grainmodel}
$$
\Phi(\P):= \bigcup_{\left\{ \substack{w\in \R^{d-1} \times\R \\
\P \cap \rm{int}(\Pi^{\downarrow}(w)) =\emptyset}\right.} \Pi^{\downarrow}(w).
$$
Remove all points of $\P$  not belonging to $\partial \Phi(\P)$.  The points in $\P$ which survive this parabolic thinning are `extremal points' in $\P$; the resulting thinned point set is denoted by ${\rm{Ext}}(\P)$.
Notice that the festoon $\partial \Phi(\P)$ is a union of inverted parabolic surfaces.

When $w \in {\rm{Ext}}(\P^{( \infty, \alpha)})$, we define the scaling limit $k$-face functional $\xi_k^{( \infty)}(w, \P^{( \infty, \alpha)}):= \xi_k^{( \infty)}(w, \P^{( \infty, \alpha)} \cup \{w \} ), \ k \in \{0, 1,...,d -1 \}$, to be the product of $(k + 1)^{-1}$ and the number of $k$-dimensional faces of the festoon $\partial( \Phi( \P^{( \infty, \alpha)}))$ which contain $w$.  Otherwise, if  $w \notin {\rm{Ext}}(\P^{( \infty, \alpha)})$, then we define
$\xi_k^{( \infty)}(w, \P^{( \infty, \alpha)})$ to be zero.

We now establish explicit asymptotics, starting with expectations for the number of $k$-dimensional faces of the two random polytopes.

\begin{theo} \label{expectasy} (expectation asymptotics for the $k$-face functional).
\noindent
For all $\alpha \in (- \infty, \infty)$ and all $k \in \{0,1,...,d-1 \}$ we have
\begin{align} \label{limex2}
\lim_{n \to \infty} \frac{ \E f_k(K_{n, \alpha}) } { d \kappa_d u_{n,\alpha}^{d - 1}} &= \lim_{\la \to \infty} \frac{ \E f_k(K_{\la,\alpha}) } { d \kappa_d u_{\la,\alpha}^{d - 1}}\nonumber\\&=\int_0^{\infty} \E \xi_k^{(\infty)}((\0,h), \P^{(\infty, \alpha) } )\mathrm{d} \nu^{(\infty, \alpha)}(h).
\end{align}
\end{theo}
For $\alpha \in (- \infty, \frac{-2} {d-1} )$, we have $\nu^{(\infty, \alpha)} = \delta_{\0}$ and  \eqref{limex2} simplifies, giving
\be \label{limex1}
\lim_{n \to \infty} \frac{ \E f_k (K_{n, \alpha}) } { n } =    \E \xi_k^{(\infty)}((\0,0), \P^{(\infty, \alpha) } ).
\ee
When $k = 0$ the right-hand side of \eqref{limex1} further simplifies and equals  $1$. Thus, in the large $n$ limit, all points remain in convex position;
this phenomenon is explained in Remark (iv) in Section 1.4.

Expectation asymptotics for the number $f_0$ of extreme points is complemented by the convergence of the set of extreme points as a point process.

\begin{theo} \label{Thm2} (convergence of extreme points).
  Let $\alpha \in (- \infty, \infty)$. As $n \to \infty$ (resp. $\la\to\infty$), under the transformation $T^{(n, \alpha)} $ (resp. $T^{(\la,\alpha)}$), the extreme points of $K_{n, \alpha}$ (resp. $K_{\la,\alpha}$) converge in distribution to $\rm{Ext}( \P^{(\infty, \alpha) } )$.
\end{theo}

For all $v \in \R^{d-1}$ and $r \in (0, \infty)$ we let ${\cal C}(B_{d-1}(v, r))$ be the space of continuous functions on $B_{d-1}(v, r)$ equipped with the supremum norm. We prove the convergence of the boundary of each random polytope, seen as a random closed set.

\begin{theo} \label{Thm3} (convergence of convex hull boundary). Fix $L \in (0, \infty)$.  For all  $\alpha \in (-\infty, \infty)$, the  rescaled boundary $T^{(\lambda, \alpha)}(\partial K_{\la,\alpha} )$ converges in distribution  $\lambda \to \infty$ to $\partial( \Phi(\P^{(\infty, \alpha) } ))$ in the space ${\cal C}(B_{d-1}(\0, L))$.
When $\alpha \in (-\infty, - \frac{2}{d-1})$, the convergence holds a.s.
\end{theo}
We now turn to the second-order description of the number of $k$-dimensional faces of random polytopes in the Poisson model. We start by introducing the proper candidate for the limiting variance.
\vskip.1cm
\begin{defn} For all $\v_1, \v_2 \in \R^{d-1} \times \R^+$ and all $\xi_k^{(\infty, \alpha)}, k \in \{0,1,...,d-1 \}$, $\alpha \in (-\infty, \infty)$  put
\begin{align} \label{SO2aa} c^{\xi_k^{(\infty, \alpha)} }(\v_1,\v_2)
& :=  \E \xi_k^{(\infty, \alpha)}(\v_1, \P^{(\infty, \alpha)}
\cup \{\v_2\} ) \xi_k^{(\infty, \alpha)}(\v_2, \P^{(\infty, \alpha)} \cup \{\v_1\} )\nonumber \\
&  \ \ \ \ \ \ \ -  \E
\xi_k^{(\infty, \alpha)}(\v_1, \P^{(\infty, \alpha)} ) \E \xi_k^{(\infty, \alpha)}(\v_2, \P^{(\infty, \alpha)} )
 \end{align}
and
\begin{align} \label{S03}
 \sigma^2(\xi_k^{(\infty, \alpha)}) & :=  \int_{0}^{\infty} \E \xi_k^{(\infty, \alpha)}((\0,h), \P^{(\infty, \alpha)} )^2 \mathrm{d}\nu^{(\infty, \alpha)}(h)  \nonumber \\
  &    \ \  \ + \int_{0}^{\infty}  \int_{\R^{d-1}} \int_{0}^{\infty} c^{\xi_k^{(\infty, \alpha)}}((\0,h),(v_1,h_1))\mathrm{d} \mu^{(\infty, \alpha)}(v_1,h_1)\mathrm{d} \nu^{(\infty,\alpha)}(h).
 \end{align}
\end{defn}
When $\alpha<-\frac{2}{d-1}$, we obtain the further simplification
$$\sigma^2(\xi_k^{(\infty, \alpha)}) =\E \xi_k^{(\infty, \alpha)}((\0,0), \P^{(\infty, \alpha)} )^2+\int_{\R^{d-1}}c^{\xi_k^{(\infty, \alpha)}}((\0,0),(v_1,0))\mathrm{d}v_1.$$
%\vskip.5cm
The next theorem states precise variance asymptotics which are specific to the Poisson model.
\begin{theo} \label{variance} (variance  asymptotics for the $k$-face functional).
For all $\alpha \in (- \infty, \infty)$ and  all $k \in \{0,1,...,d-1 \}$ we have
\be \label{lim1}
\lim_{\la \to \infty} \frac{ \Var f_k(K_{\la,\alpha}) } {d \kappa_d  u_{\la,\alpha}^{d - 1}} = \sigma^2(\xi_k^{(\infty, \alpha)}) \in (0, \infty).
\ee
\end{theo}
When $\alpha \in (-\infty, \frac{-2} {d-1} )$, we have by  Table \ref{tab:summary} that  $d \kappa_d \cdot u_{\la,\alpha}^{d - 1} = \la$ and thus
\be \label{lim1aa}
\lim_{\la \to \infty} \frac{ \Var f_k(K_{\la,\alpha}) } {\la } =  \sigma^2(\xi_k^{(\infty, \alpha)}).
\ee
When $k= 0$ we obtain the further simplification
\be \label{lim1a}
\lim_{\la \to \infty} \frac{ \Var f_0(K_{\la,\alpha}) } {\la } = 1. \ee

Finally, we obtain rates of asymptotic normality for $k$-face functional.

\vskip.1cm

\begin{theo} \label{clt} (rates of normal convergence for $k$-face functional).  For all $\alpha \in ( - \infty, \infty)$
and all $k \in \{0,1,...,d-1 \}$ we have
$$
\sup_{t} \left| \PP( \frac{  f_k(K_{\la,\alpha}) -   \E f_k(K_{\la,\alpha})    }  { \sqrt{ \Var f_k(K_{\la,\alpha}) } }  \leq t ) - \PP( {\cal N} \leq t) \right|
= O\left(  \frac{ (\log \la)^{3d + 1}}  {   \sqrt{\Var f_k(K_{\la,\alpha}) }}\right),
$$
where ${\cal N}$ denotes a mean zero normal random variable with variance one.
\end{theo}

\vskip.2cm
The reader may wonder about the connection between this paper and the authors' earlier papers.  Here we show that some, but not all of the techniques of \cite{CSY, CY2, SY} may be used to establish the limit theory for a new random polytope model. In this new model, the input is not hosted by a smooth mother body as in \cite{CSY, CY2, SY}, but is instead hosted by a union of balls of radius $n^{\alpha}$ with centers randomly distributed on the unit sphere. Unlike previously studied models, this new model exhibits phase transitions, multiple scaling regimes, and is asymptotically described by not just one, but five limit point processes.

A further limitation of the papers \cite{CSY, CY2, SY} is that they do not cover the case where
 the limit point process is hosted by a subset of $\R^{d - 1} \times [0, \infty)$, including e.g. $\R^{d - 1} \times [0,2]$ and $\R^{d - 1} \times \{\0 \}$, which arise when $\alpha = \frac{-2} {d-1}$ and $\alpha \in (- \infty, \frac{-2} {d-1})$, respectively.  When the limit point process is hosted by $\R^{d - 1} \times \{\0 \}$, i.e., when  $\alpha \in (- \infty, \frac{-2} {d-1})$, then we require new techniques leading to strengthened results concerning uniform convergence of one and two point correlation functions.
 When  $\alpha \in (\frac{-2} {d-1} , \frac{2}{d + 1} ) \cup (\frac{2}{d + 1} , \infty)$ the limit point processes ${\cal P}^{(\infty, \alpha)}$ coincide with those in \cite{CSY}, even though the input is not of the same form.
Additional new aspects of this paper include de-Poissonization and proof of the strict positivity of limiting variances.
\vskip.3cm

\noindent{\bf 1.4.  Remarks.}
(i) {\em Comparison with the literature}. The scaling exponents generating the scaling regimes
 of Theorem \ref{expectasy}  are consistent with the rates of growth for the $k$-face functional in Theorem 7 of \cite{DGGT}.  That theorem gives rates of growth whereas we find explicit asymptotics as well, including the case of critical $\alpha$ values.  It is a curious fact that, at least in the large $n$ limit, when $\alpha \in (-\infty, \frac{-2} {d-1})$
 the average number of $k$-faces  equals the worst-case number, namely $n$,  and then, as $\alpha$ increases up to $0$, the average number decreases down to the minimum number, which is of order $\Theta(n^{(d-1)/2d})$.   Thereafter, for increasing positive $\alpha$, the order of the average number of $k$-faces increases up to $\Theta(n^{(d-1)/(d + 1)})$, which is precisely the growth rate achieved in the average-case analysis for the standard model consisting of Poisson input on the unit ball \cite{CSY, CY, CY2}.
Indeed, for all $\alpha \in ( \frac{2} {d + 1}, \infty)$  both the expectation and variance asymptotics for the $k$-face functional on perturbed input $\tilde{\P}_{\la, \alpha}$  coincide with their counterparts involving  Poisson input on the unit ball \cite{CSY, CY, CY2}.

\vskip.1cm

\noindent(ii) {\em The case $\alpha \in (\frac{2} {d + 1}, \infty)$}.
In this regime, it is more appropriate to view the model by exchanging the roles of the initial input and the perturbation, i.e., to consider the sequence $\{e_i\}_{i=1}^n$ as the initial input and the sequence $\{X_i\}_{i=1}^n$ as the perturbation. Classical results tell us that  the convex hull of $n^{-\alpha} \{e_i \}_{i = 1}^n$ is with high probability contained in an annulus of width $\Theta(n^{-2/(d + 1)})$ and the number of extreme points is $\Theta(n^{(d - 1)/(d + 1)} )$.  Consequently, with high probability the width of the annulus holding the convex hull of $\{e_i \}_{i = 1}^n$ is $\Theta(n^{\alpha} n^{-2/(d + 1)})$, which exceeds the diameter of $\{X_i\}_{i=1}^n$ when $\alpha \in ( \frac{2} {d + 1}, \infty)$.  Thus, `perturbing' each $e_i, 1 \leq i \leq n,$ by $X_i$ modifies this annulus by a negligible amount, as $|X_i| = 1$. One might expect that it also modifies the cardinality of extreme points by a negligible amount.

 \vskip.1cm

\noindent(iii)  {\em The case $\alpha = 0$}.  As mentioned, Theorem \ref{expectasy} shows that the expected number of extreme points is minimized when $\alpha = 0$.  In this case the order of the expected number of extreme points asymptotically behaves like $n^{(d-1)/2d}$. The case $\alpha=0$ is precisely the situation when the perturbed region contains the origin on its boundary for all $n$. The probability that a perturbed point is close to the origin and in particular interior to the convex hull is then maximal, which intuitively explains why the number of extreme points should be minimal.
If one allows perturbations in $B_d(\0, \rho)$, $\rho$ a constant, then straightforward modifications of our approach show that the expectation asymptotics in Theorem \ref{expectasy} involve an additional factor of $((1 + \rho)^d \times 2^{-d} \times \rho^{-(d + 1)/2})^{(d-1)/2d}$, which is minimized when $\rho = \frac{d + 1} {d - 1}$. {\em In other words, as $n \to \infty$, the expected value of the $k$-face functional on perturbed input achieves its minimum when the perturbations are confined to a ball of radius $\rho = \frac{d + 1} {d - 1} $.}

\vskip.1cm

\noindent(iv) {\em The case $\alpha \in (-\infty, \frac{-2} {d - 1})$}.   The minimal interpoint distance between the non-perturbed points  $\X_n$ is with high probability $\Omega(n^{-2/(d - 1)})$, as shown in \cite{JJ}.   When the points undergo perturbations with magnitude $o(n^{-2/(d - 1)})$, as is the case when $\alpha \in (-\infty, \frac{-2} {d - 1} )$, then simple geometry shows that for large $n$, all the points remain in extreme position. This is the content of \eqref{limex1}, and, in the case of Poisson input, of \eqref{lim1a}.

\vskip.1cm

\noindent(v) {\em  Geometric significance of $\beta$.}  The  scaling exponents $\beta$ in Figure \ref{fig:graphbeta}
 describe two geometric properties. First, up to a constant,  $n^{\alpha} n^{-2 \beta}$ represents the expected distance between a point of $\partial K_{n, \alpha}$ and the boundary of $B_d(\0, 1 + n^{\alpha})$.
Second, the average number of extreme points is of the order $n^{\beta (d - 1)}$.

\vskip.1cm

\noindent(vi) {\em a.s. convergence.}   Apart from Theorem \ref{Thm3}, we have not sought a.s. convergence results.  Yet our proof techniques yield a.s. convergence in some instances.  For example, as seen in the proof of
Proposition \ref{expectation-1} below, for each  $\alpha \in (- \infty, \frac{-2} {d-1} )$ there is a coupling of $\P_{\la}, \la \geq 1,$ and $\P^{(\infty, \alpha)}$  such that a.s. $T^{(\la, \alpha)} (\tilde{\P}_{\la, \alpha}) \to \P^{(\infty, \alpha) }.$  Thus Proposition \ref{prop1} may be strengthened for $\alpha \in (-\infty, \frac{-2} {d-1})$.
Likewise, this coupling shows that the extreme points of $K_{\lambda,\alpha}$ converge a.s. to $\rm{Ext}( \P^{(\infty, \alpha)})$.

\vskip.1cm

\noindent(vii) {\em Extensions.}
We expect that versions of Theorems \ref{variance} and \ref{clt} hold for binomial input.  For $\alpha>\frac{2}{d+1}$, it is likely that  the limit \eqref{lim1} may be de-Poissonized, as it would be similar to what has already been done for uniform input in a smooth convex body as in \cite{CY}.   However, for $\alpha<-\frac{2}{d-1}$, we anticipate that  the limiting variance for the binomial
model to be zero, hence smaller than the limiting variance at \eqref{lim1aa}.
 We expect that our main results also hold for defect and intrinsic volumes of the convex hull $K_{\la,\alpha}$.  We leave these questions for future research.

\vskip.1cm

The paper is structured as follows: Sections \ref{sec:scalinglimits}, \ref{sec:stab} and \ref{sec:proofs} are devoted to the proofs of Proposition \ref{prop1}, several intermediary results on stabilization in the rescaled space and the proofs of Theorems \ref{expectasy}--\ref{clt} respectively.

\section{Finite-size scaling and scaling limits}\label{sec:scalinglimits}

\allco

The aim of this section is to prove Proposition \ref{prop1}, i.e., to establish the convergence of the transformed perturbed binomial and Poisson point processes to the limiting Poisson point process $\P^{(\infty,\alpha)}$. We start by studying the intensity measure of these processes in Subsection 2.1 before moving on to the asymptotics in Subsection 2.2.
Denote the closure of the injectivity region
%$\{ \tilde{v} \in T(\S^{d-1}), | \tilde{v} | \leq \pi \} $
of the exponential map by $B_{d-1}(\pi)$.
In view of the mapping \eqref{map2} we put
$$
v := u_{\la,\alpha} \exp_{d-1}^{-1} \frac{x} {|x|}, \ \ \ h := u_{\la,\alpha}^2(1 - \frac{|x|}{ 1 + \la^{\alpha} } ).
$$
Note that $v$ ranges over $u_{\la,\alpha} B_{d-1}(\pi)$. A consequence of the upcoming Lemma \ref{lem:propfn} is that $h$ ranges over $[0, h^{(\la,\alpha)} ]$, where
\be \label{udef}
h^{(\la,\alpha)}  := \left\{ \begin{array}{ll} u_{\la,\alpha}^2 &  ~  \alpha \in [ -\frac{2} {d - 1} , \infty)\\
\frac{ 2\la^{\alpha} u_{\la,\alpha}^2} { 1 + \la^{\alpha} }
&  ~ \alpha \in (- \infty, -\frac{2}{d-1} ). \\
\end{array}
\right.
\ee
\vskip.1cm
\noindent{\bf 2.1. Measures of caps.}
To derive properties of the intensity measure of the Poisson point process $\P^{(\la, \alpha)}$,
 we shall need some measure theoretic estimates for the normalized Hausdorff measure of caps on $\S^{d-1}$.
Here and elsewhere $c, c_1, c_2...$ denote generic finite positive constants whose value may change at each occurrence and which may depend on $d$ and $\alpha$, unless stated otherwise.

\begin{lemm}\label{lem:propfn}
The following claims hold for all $\alpha \in (-\infty, \infty)$. \\
\noindent (i) Let $x\in\R^d$ be at distance $(1+\la^{\alpha})(1-u_{\la,\alpha}^{-2}h)$ from the origin. The normalized Hausdorff measure of the intersection $B_d(x,\la^{\alpha})\cap \S^{d-1}$ equals
$$\frac{1}{d\kappa_d} \sigma_{d-1}(B_d(x,\la^{\alpha})\cap \S^{d-1})=\stwo(g_{\la, \alpha} (h))$$
where the function $\stwo$ is defined at \eqref{eq:defs2} and where
\begin{equation}
  \label{eq:defgnh}
g_{\la, \alpha} (h):=\left(\frac{\la^{\alpha}u_{\la,\alpha}^{-2}h-\frac{1}{2}(1+\la^{\alpha})u_{\la,\alpha}^{-4}h^2}{1-u_{\la,\alpha}^{-2}h}\wedge 2\right)\vee 0.
\end{equation}
In particular,
 \begin{equation*}
 \stwo(g_{\la, \alpha} (h))=1 \Longleftrightarrow g_{\la, \alpha} (h)=2 \Longleftrightarrow  (\alpha>0) \mbox{ and } (h \in  [\frac{2u_{\la,\alpha}^{2}}{1+\la^{\alpha}}, u_{\la,\alpha}^2))
 \end{equation*}
 and where
  \begin{equation*}
 \stwo(g_{\la, \alpha} (h))=0 \Longleftrightarrow g_{\la, \alpha} (h)=0 \Longleftrightarrow  (\alpha<0) \mbox{ and } (h \in [ h^{(\la,\alpha)}, u_{\la,\alpha}^2)).
  \end{equation*}
%$g_{\la, \alpha} (h)=2$ iff $\alpha>0$ and $h \in  [\frac{2u_{n,\alpha}^{2}{1+n^{\alpha}}, u_{n,\alpha}^2)$ in which case $\stwo(g_{\la, \alpha} (h))=d\kappa_d$ and $g_{\la, \alpha} (h)=0$ iff $\alpha<0$ and $h \in [ \un, u_{n,\alpha}^2)$ in which case $\stwo(g_{\la, \alpha} (h))=0$. \\~\\
%\Comment{ $h \in ( \un, u_{n,\alpha}^2)$;  is it ok? Details added. OK? PC}
\noindent (ii)
If $g_{\la, \alpha} (h)\to 0$ when $\la\to \infty$, then
$$\stwo(g_{\la, \alpha} (h))\sim \frac{\kappa_{d-1}}{d\kappa_d}(2 g_{\la, \alpha} (h))^{\frac{d-1}{2}}.$$
\noindent (iii) There exist constants $c_1$ and $c_2$ such that for  any $h\in (0, u_{\la,\alpha}^2)$,
\begin{align}
  \label{eq:boundsf}
&c_1 \left((\la^{\alpha}u_{\la,\alpha}^{-2}h((1-\frac{1+\la^{-\alpha}}{2}u_{\la,\alpha}^{-2}h)\vee 0))\wedge 1\right)^{\frac{d-1}{2}}
\nonumber\\&\hspace*{2cm}
\le
\stwo(g_{\la, \alpha} (h))
\le
\left( c_2 (\la^{\frac{\alpha}{2}}u_{\la,\alpha}^{-1})^{d-1}\frac{h^{\frac{d-1}{2}}}{(1-u_{\la,\alpha}^{-2}h)^{\frac{d-1}{2}}}\right)\wedge 1.
\end{align}
\end{lemm}
\noindent{\em Proof}.
(i) We start by treating the two particular cases. The normalized Hausdorff measure of $B_d(x,\la^{\alpha})\cap \S^{d-1}$ is equal to $1$ if and only if $B_d(x,\la^{\alpha})$ contains the entirety of $\S^{d-1}$. This is equivalent to $|x| \leq \la^{\alpha} - 1$, i.e.,
$h$ is such that
 $(1+\la^{\alpha})u_{\la,\alpha}^{-2}h  \in [2, \infty)$. Since $h\in (0,u_{\la,\alpha}^{2})$, this implies that $\alpha>0$.
We finally remark that the inequality $(1+\la^{\alpha})u_{\la,\alpha}^{-2}h  \geq 2$
holds  for $\la$ large enough as soon as $\alpha-2\beta > 0$, which is the case when $\alpha \in ( \frac{2}{d+1}, \infty)$.

Similarly, the normalized Hausdorff measure of $B_d(x,\la^{\alpha})\cap \S^{d-1}$ is equal to $0$ if and only if $B_d(x,\la^{\alpha})$ is included in the interior of  $B_d({\bf 0},1)$. This is equivalent to
$|x| < 1 - \la^{\alpha}$, i.e.,
$(1+\la^{\alpha})u_{\la,\alpha}^{-2}h>2\la^{\alpha}$.
Since $h\in (0,u_{\la,\alpha}^2)$, this implies that $\alpha<0$. We finally remark that the inequality $(1+\la^{\alpha})u_{\la,\alpha}^{-2}h>2\la^{\alpha}$ occurs for $\la$ large enough as soon as $\alpha+2\beta < 0$, which is the case if $\alpha \in (-\infty, -\frac{2}{d-1})$.

Assume now that we are in the general case when $B_d(x,\la^{\alpha})$ has a non-trivial intersection with $\S^{d-1}$. For $\la$ large enough, this happens only when  $\alpha \in [ \frac{-2} {d-1}, \frac{2} {d + 1}].$
The normalized Hausdorff measure of $B_d(x,\la^{\alpha})\cap \S^{d-1}$ equals $\stwo(1-\cos(\phi_\la(h)))$ where
$\phi_\la(h)$ is the angle between the first and second edges of a triangle with edge lengths $1$, $(1+\la^{\alpha})(1-u_{\la,\alpha}^{-2}h)$ and $\la^{\alpha}$. The law of cosines yields
\begin{align*}
1-\cos(\phi_\la(h))&=1-\frac{1+(1+\la^{\alpha})^2(1-u_{\la,\alpha}^{-2}h)^2-\la^{2\alpha}}{2(1+\la^{\alpha})(1-u_{\la,\alpha}^{-2}h)}\nonumber\\
&=\frac{\la^{\alpha}u_{\la,\alpha}^{-2}h-\frac{1}{2}(1+\la^{\alpha})u_{\la,\alpha}^{-4}h^2}{1-u_{\la,\alpha}^{-2}h} %\label{eq:formulacosphi}
\end{align*}
where the second equality follows by expanding the numerator and re-writing it as
$$
2(1 + \la^{\alpha} )(1-u_{\la,\alpha}^{-2}h) + u_{\la,\alpha}^{-2} h(-2\la^{\alpha}(1+\la^{\alpha}) + u_{\la,\alpha}^{-2} h(1 + \la^{\alpha})^2).
$$
Now $g_{\la, \alpha} (h)= \left((1-\cos(\phi_\la(h)) \wedge 2\right)\vee 0$ and noting that we have already treated the cases $g_{\la, \alpha}  = 2$ and $g_{\la, \alpha}  = 0$, we obtain (i).

\vskip.1cm

\noindent(ii) This estimate is deduced from (i) using the equivalences $\arccos (1-u)\underset{u\to 0}{\sim} \sqrt{2u}$ and $\int_0^{\phi}\sin^{d-2}\theta\mathrm{d}\theta\underset{\phi\to 0}{\sim} \frac{1}{d-1}\phi^{d-1}$.
\vskip.2cm
\noindent (iii)
We only need to prove \eqref{eq:boundsf} in the case when $\stwo(g_{\la, \alpha} (h)) \in (0, 1)$. Recall that  $\sqrt{2u}\le \arccos(1-u)\le 2\sqrt{2u}$ for every $u\in [0,2]$. This implies that
$$
\sqrt{(2\la^{\alpha}u_{\la,\alpha}^{-2}h-(1+\la^{\alpha})u_{\la,\alpha}^{-4}h^2)\vee 0}\wedge 2 \le \arccos(1-g_{\la, \alpha} (h)) \le \frac{2\sqrt{2h} \la^{\frac{\alpha}{2}}u_{\la,\alpha}^{-1}}{\sqrt{1-u_{\la,\alpha}^{-2}h}}.
$$
Moreover, we have the following upper and lower  bounds %for the integral of $\sin^{d-2}(\theta)$:
$$
\int_0^{\phi}\sin^{d-2}\theta\mathrm{d}\theta\le \int_0^{\phi}\theta^{d-2}\mathrm{d}\theta=\frac{1}{d-1}\phi^{d-1}
$$
and
$$
\int_0^{\phi}\sin^{d-2}\theta\mathrm{d}\theta\ge \int_0^{\phi/2}\sin^{d-2}\theta\mathrm{d}\theta\ge \int_0^{\phi/2}\left(\frac{2 \theta }{\pi} \right)^{d-2}\theta\mathrm{d}\theta=\frac{\phi^{d-1} }{2\pi^{d-2}(d-1)}.
$$
Combining \eqref{eq:defs2} and the last three displays yields \eqref{eq:boundsf}. \qed

~\\

\noindent{\bf 2.2. The intensity measure of perturbed points.} Let $\mu^{(\la,\alpha)}$ be the intensity measure of the rescaled Poisson point process $\P^{(\la,\alpha)}$.  Note that $\mu^{(\la,\alpha)}$ is defined on $u_{\la,\alpha}B_{d-1}(\pi)\times[0,h^{(\la,\alpha)}]$.
We  establish a formula for the density of $\mu^{(\la,\alpha)}$, an asymptotic equivalence, and upper and lower bounds.

\begin{lemm} \label{lem:proppn}
Let $\alpha \in (-\infty, \infty)$.
\\
\noindent (i)  The measure $\mu^{(\la,\alpha)}$  has a density given by
 \begin{align}  \label{eq:formulapnf}
\frac{\sin^{d-2}(u_{\la,\alpha}^{-1}|v|)}{|u_{\la,\alpha}^{-1}v|^{d-2}}\varphi^{(\lambda,\alpha)}(h)
\end{align}
where for any $h\in [0,h^{(\lambda,\alpha)}]$,
% where $\mu_v^{(\lambda,\alpha)}$ is a measure on $u_{\la,\alpha}B_{d-1}(\pi)$ and and $\varphi^{(\lambda,\alpha)}(h)$ is a function on $[0,h^{(\lambda,\alpha)}]$ respectively defined as
 \begin{align} \label{defph}
\varphi^{(\lambda,\alpha)}(h)& :=\frac{\lambda(1+\lambda^{\alpha})^{d}u_{\la,\alpha}^{-(d+1)}}{\kappa_d \lambda^{d\alpha}}\stwo(g_{\la, \alpha} (h))(1-u_{\la,\alpha}^{-2}h)^{d-1}.
\end{align}
In particular,
\begin{align}
  \label{eq:totalmassmulambdaalpha}
\lambda&%=\int_0^{h^{(\lambda,\alpha)}}\int_{u_{\la,\alpha}B_{d-1}(\pi)} \mathrm{d}\mu_v^{(\lambda,\alpha)}(v)\varphi^{(\lambda,\alpha)}(h)\mathrm{d}h
=d\kappa_d u_{\la,\alpha}^{d-1}\int_0^{h^{(\lambda,\alpha)}}\varphi^{(\lambda,\alpha)}(h)\mathrm{d}h.
\end{align}
\noindent (ii) If $g_{\la, \alpha} (h)\to 0$ when $\la\to \infty$, then
\begin{equation}
  \label{eq:equivintermpn}
\varphi^{(\la,\alpha)}(h)\sim \frac{(1+\la^{\alpha})^d\kappa_{d-1}}{d\kappa_d^2}\la^{1-d\alpha}u_{\la,\alpha}^{-(d+1)}(2g_{\la, \alpha} (h))^{\frac{d-1}{2}}.
\end{equation}
\noindent (iii)  Put $e(\la):= u_{\la,\alpha}^{-(d+1)}\la^{(d\alpha\vee 0) +1-d\alpha}$.
There exist constants $c_1, c_2 \in (0, \infty)$ such that for any $0<h_1<h_2< \un$,
\begin{align}
 & c_1 e(\la)[(\la^{\frac{\alpha}{2}}u_{\la,\alpha}^{-1})^{d-1}h^{\frac{d-1}{2}}((1-\frac{1+\la^{-(\alpha\wedge 0)}}{2}u_{\la,\alpha}^{-2}h)\vee 0)^{\frac{3(d-1)}{2}}\wedge (1-u_{\la,\alpha}^{-2}h)^{d-1}]
%c_1  n^{d\alpha\vee 0+1-(d+1)\frac{\alpha}{2}}u_{\la,\alpha}^{-2d}h_1^{\frac{d-1}{2}}\left((1-\frac{1+n^{-(\alpha\wedge 0)}}{2}u_{\la,\alpha}^{-2}h_2)\vee 0\right)^{\frac32(d-1)}
%(1-u_{\la,\alpha}^{-2}h_2)^{d-1}\left((\min_{h\in [h_1,h_2]}(h-(1+\la^{\alpha})n^{-\alpha}u_{\la,\alpha}^{-2}\frac{h^2}{2}))\vee 0\right)^{\frac{d-1}{2}}
\nonumber  \\&  \nonumber
\le \varphi^{(\la,\alpha)}(h) \nonumber \\&
\le c_2 e(\la) \left(\left((\la^{\frac{\alpha}{2}}u_{\la,\alpha}^{-1})^{d-1}
h^{\frac{d-1}{2}}\right)\wedge 1\right). \ \  \label{eq:boundspn}
\end{align}
\end{lemm}

\vskip.2cm
\noindent{\em Proof}.
(i) Let us first determine the distribution of the rescaled perturbed random point $T^{(\la,\alpha)}(\tilde{X}_1)$. For $v \in u_{\la,\alpha} B_{d-1}(\pi)$, $0 < h_1 <h_2 < h^{(\la,\alpha)}$
and $r \in (0, u_{\la,\alpha}\pi)$ we put
 $$
 A_\la(v,r,h_1,h_2):=[T^{(\la, \alpha)} ]^{-1}(B_{d-1}(v,r) \times [h_1, h_2]).
 $$
Conditional on $X_1$, the distribution of  $\tilde{X}_1$  is uniform inside $B_d(X_1,\la^{\alpha})$. Consequently, the probability that $\tilde{X}_1$ falls inside the set  $A_\la(v,r,h_1,h_2)$ is the ratio of the expected volume of $B_d(X_1,\la^{\alpha})\cap A_\la(v,r,h_1,h_2)$ to that  of $B_d(X_1,\la^{\alpha})$.  We let $\Vol_d(A)$ denote the Lebesgue measure of a full $d$-dimensional subset of $\R^d$.
 Fubini's theorem implies
\begin{align*}
& \mathbb{P}(T^{(\la,\alpha)}(\tilde{X}_1)\in B_{d-1}(v,r)\times [h_1,h_2]) \\
&= \E \left(  \frac{ \Vol_d(B_d(X_1,\la^{\alpha})\cap A_\la(v,r,h_1,h_2))}{\Vol_d(B_d(X_1,\la^{\alpha}))} \right) \nonumber\\
&=\frac{1}{\kappa_d \la^{d\alpha}}\int_{A_\la(v,r,h_1,h_2)} \mathbb{P}  (x\in B_d(X_1,\la^{\alpha}))\mathrm{d}x\nonumber\\
&=\frac{1}{\kappa_d \la^{d\alpha}}\int_{A_\la(v,r,h_1,h_2)} \mathbb{P} (X_1 \in B_d(x,\la^{\alpha}))\mathrm{d}x.
\end{align*}
Recalling that $X_1$ is uniformly distributed on $\S^{d-1}$ we get
\begin{align}
& \mathbb{P}(T^{(\la,\alpha)}(\tilde{X}_1)\in B_{d-1}(v,r)\times [h_1,h_2]) \nonumber\\
&= \frac{1}{\kappa_d \la^{d\alpha}}\int_{A_\la(v,r,h_1,h_2)}\frac{\sigma_{d-1}(B_d(x,\la^{\alpha})\cap \S^{d-1})}{d\kappa_d}\mathrm{d}x\nonumber\\
&=\frac{\sigma_{d-1}(\exp_{d-1}(u_{\la,\alpha}^{-1}B_{d-1}(v,r)))}{\kappa_d \la^{d\alpha}} \nonumber \\
& \ \ \ \times \int_{h_1}^{h_2}\stwo(g_{\la, \alpha} (h))(1+\la^{\alpha})^{d}(1-u_{\la,\alpha}^{-2}h)^{d-1}u_{\la,\alpha}^{-2}\mathrm{d}h, \label{pnformula}
\end{align}
where  $|x| = (1 + \la^{\alpha})(1 - u_{\la,\alpha}^{-2} h)$ and $\mathrm{d}x = |x|^{d-1} (1 + \la^{\alpha} ) u_{\la,\alpha}^{-2} \mathrm{d}h \mathrm{d} \sigma_{d-1}$.

We may derive from \eqref{pnformula} the density of $T^{(\la,\alpha)}(\tilde{X}_1)$ by recalling
 the Jacobian of the exponential map $\exp_{d-1}(u_{\la, \alpha}^{-1} \cdot )$, which is precisely the product of $u_{\la, \alpha}^{-(d-1)}$ and the first factor in \eqref{eq:formulapnf}; c.f.  \cite[display (2.17) with $u_{\la,\alpha}^{-1}$ in place of $\la^{-\beta}$]{CSY}.
 Since the intensity measure $\mu^{(\la,\alpha)}$ equals the distribution of $T^{(\la,\alpha)}(\tilde{X}_1)$ multiplied by $\la$, we obtain  \eqref{eq:formulapnf}. Remembering that the total mass of $\mu^{(\la,\alpha)}$ is $\la$, we obtain \eqref{eq:totalmassmulambdaalpha}.

\vskip.2cm

\noindent (ii) The equivalence \eqref{eq:equivintermpn} follows from the definition of $\varphi^{(\lambda,\alpha)}$ at \eqref{defph}  and  Lemma \ref{lem:propfn}(ii).

\vskip.2cm
\noindent (iii)
We first show the upper bound in \eqref{eq:boundspn}.  Combining (i) with  the right-hand side of \eqref{eq:boundsf} gives
\begin{align*}
 \varphi^{(\la,\alpha)}(h)
&\le c u_{\la,\alpha}^{-(d+1)}\la^{1-d\alpha}(1+\la^{\alpha})^d\left[\left((\la^{\frac{\alpha}{2}}u_{\la,\alpha}^{-1})^{d-1}h^{\frac{d-1}{2}}(1-u_{\la,\alpha}^{-2}h)^{\frac{d-1}{2}}\right)\wedge 1\right]
\end{align*}
which yields the right-hand side of
\eqref{eq:boundspn}. We show now the left-hand side of \eqref{eq:boundspn}.
Combining the lower bound for  $\stwo(g_{\la, \alpha} (h))$ at \eqref{eq:boundsf} with  \eqref{eq:formulapnf}
gives
\begin{align*}
& \varphi^{(\la,\alpha)}(h)
\ge c e(\la)\left[\left((\la^{\alpha}u_{\la,\alpha}^{-2}h(1-\frac{1+\la^{-\alpha}}{2}u_{\la,\alpha}^{-2}h)\vee 0)^{\frac{d-1}{2}}(1-u_{\la,\alpha}^{-2}h)^{d-1}\right)\right.\\&\hspace*{9cm}\left.\wedge (1-u_{\la,\alpha}^{-2}h_2)^{d-1}\right]
\end{align*}
which yields  the left-hand side of \eqref{eq:boundspn} and concludes the proof of Lemma \ref{lem:proppn}.
\qed

\vskip.2cm
%\noindent{\bf 2.2.  The limit intensity measure of perturbed points.}
The following lemma figures prominently in the proof of  Proposition \ref{prop1} and it provides
 useful lower and upper bounds on the density of $\mu^{(\la,\alpha)}$. When $\alpha\in [-\frac{2}{d-1},\infty)$, we denote by $\varphi^{(\infty,\alpha)}(h)$ the density of $\nu^{(\infty,\alpha)}$.

\vskip.2cm

\begin{lemm} \label{densitybounds} (i) For all $\alpha \in (-\infty, \infty)$ we have as $\la \to \infty$ the convergence of measures
\begin{align}  \label{generalbounda}
\mu^{(\la,\alpha)} \to\mu^{(\infty, \alpha)}
\end{align}
where $\mu^{(\infty,\alpha)}$ is defined at \eqref{eq:deflimitingintensity}.

\noindent (ii) For all $\alpha\in [-\frac{2}{d-1},\infty)$ and $h\in [0,h^{(\lambda,\alpha)}]$,  we have as $\la \to \infty$
\begin{equation}
  \label{eq:convphi}
\varphi^{(\la,\alpha)}(h)\to \varphi^{(\infty,\alpha)}(h).
\end{equation}
Moreover,
\begin{equation}
  \label{eq:upperboundphi}
\varphi^{(\la,\alpha)}(h)\le \left\{\begin{array}{ll} c & ~ \alpha \in [\frac{2} {d + 1}, \infty)\\
c h^{\frac{d-1}{2}} & ~\alpha \in[-\frac{2}{d-1},\frac{2}{d+1})
\end{array}\right.
\end{equation}
and for $v\in \frac{u_{\la,\alpha}}{2}B_{d-1}(\pi)$ and $h\in [\frac{1}{4}\kappa_d^{ -2/ (d+1) },\frac{1}{2}\la^{\alpha\wedge 0}u_{\la,\alpha}^2]$,
\begin{equation}
  \label{eq:lowerboundmu}
\frac{\mathrm{d}\mu^{(\la,\alpha)}(v,h)}{\mathrm{d}v\mathrm{d}h}\ge c.
\end{equation}
\noindent (iii) For all $\alpha \in (-\infty, \frac{-2} {d - 1})$,
\begin{equation}
  \label{eq:lowerboundmualphasmall}
c_1\le \frac{\int_{h=0}^{h^{(\lambda,\alpha)}}\mathrm{d}\mu^{(\la,\alpha)}(v,h)}{\mathrm{d}v}\le c_2.
\end{equation}
 \end{lemm}

\vskip.2cm

\noindent{\em Proof.}
We establish Lemma \ref{densitybounds} on a case-by-case basis according to the regimes for the parameter $\alpha$ described in  Table \ref{tab:summary}. We make some preliminary observations used in every case: The first term $$\frac{\sin^{d-2}(u_{\la,\alpha}^{-1}|v|)}{|u_{\la,\alpha}^{-1}v|^{d-2}}$$ in the density of $\mu^{(\la,\alpha)}$ given at \eqref{eq:formulapnf} converges to $1$ for all $v\in \R^{d-1}$; it is upper-bounded by $1$ for all $v\in u_\la B_{d-1}(\pi)$ and lower-bounded by a positive constant as soon as $u_n^{-1}|v|\le \frac{1}{2}$. This implies in particular that for $\alpha\in [-\frac{2}{d-1},\infty)$,
the convergence \eqref{generalbounda} is a direct consequence of \eqref{eq:convphi}. For $\alpha\in (-\infty,-\frac{2}{d-1})$, (i) and (iii) amount to showing that the total mass of $\varphi^{(\la,\alpha)}(h)$ converges to $1$ and is lower and upper bounded.

\vskip.2cm

\noindent {\it  (a) The regime  $\alpha \in (\frac{2} {d+1}, \infty)$.}
    Recall $\beta = \frac{1} {d + 1}$ and thus  for $\la$ large enough, we have $(1+\la^{\alpha})u_{\la,\alpha}^{-2}h \in (2, \infty)$, whence  $\stwo(g_{\la, \alpha} (h))$ tends to $1$ as $\la \to \infty$. Inserting this into \eqref{eq:formulapnf}, we get that  $\varphi^{(\la,\alpha)}(h)$ converges  to $1$ for all $h>0$, which shows \eqref{generalbounda}.

We now bound $\varphi^{(\la,\alpha)}(h)$: noticing that $e(\la)$ is a constant, we obtain \eqref{eq:upperboundphi} by considering the second term in the minimum on the right-hand side of \eqref{eq:boundspn}.

For the lower bound, we fix $h\in [\frac{1}{4}\kappa_d^{-2/(d+1)},\frac{1}{2}u_{\la,\alpha}^2]$ and consider the left-hand side of \eqref{eq:boundspn}. The first term in the minimum therein goes to $\infty$ because the quantity $\la^{\alpha/2 }u_{\la,\alpha}^{-1}$ goes to $\infty$ while the quantity  $h^{(d-1)/2}  (1-u_{\la,\alpha}^{-2}h)^{3(d-1)/2}$ is bounded from below by a constant. Consequently, that minimum is equal to the second term $(1-u_{\la,\alpha}^{-2}h)^{d-1}$ which is bounded by a positive constant. In conclusion, $\varphi^{(\la,\alpha)}(h)$ is lower-bounded by a constant and \eqref{eq:lowerboundmu} follows.
\vskip.2cm

\noindent {\it (b) The regime  $\alpha= \frac{2} {d+1}$.}
We treat the height domains $h \in [0,2\kappa_d^{-2/(d+1)})$ and $h \in [2\kappa_d^{-2/(d+1)}, \infty)$ separately.  Notice that $u_{\la,\alpha}^2 =
\kappa_d^{-2/(d+1)} \la^{\alpha}.$
For the latter domain we obtain
$$
h > \frac{2 \kappa_d^{-\frac{2}{d+1}} } { 1 + \la^{- \frac{2}{d+1} } } = \frac{2u_{\la,\alpha}^2 } { 1+\la^{\alpha} },
$$
which yields $\stwo(g_{\la, \alpha} (h))=1$ by Lemma \ref{lem:propfn} and which yields the convergence of $\varphi^{(\la,\alpha)}(h)$ to $1$.

When $h \in [0,2\kappa_d^{-2/(d+1)})$, we note that \eqref{eq:defgnh} shows that $g_{\la, \alpha} (h)\sim \kappa_d^{2/(d+1)}h$ when $\la\to\infty$, which implies in turn that $\stwo(g_{\la, \alpha} (h))$ and $\varphi^{(\la,\alpha)}(h)$ are converging to $\stwo(\kappa_d^{2/(d+1)} h)$.
In conclusion, we obtain $\varphi^{(\la,\alpha)}(h)\to \stwo(\kappa_d^{\frac{2}{d+1}}h)$,
which shows \eqref{generalbounda}.

We obtain  \eqref{eq:upperboundphi} by taking the second part of the minimum in the right-hand side of \eqref{eq:boundspn} while we derive \eqref{eq:lowerboundmu} by noticing that both parts of the minimum in the lower bound are bounded below for $h\in [\frac{1}{4}\kappa_d^{-2/(d+1)},\frac{1}{2}u_{\la,\alpha}^{2} ]$.
\vskip.2cm

\noindent {\it (c) The regime \ $\alpha\in (- \frac{2} { d - 1},\frac{2}{d+1})$.} This regimes breaks up into two sub-regimes.\\
$\cdot$ The case $\alpha\in [0,\frac{2}{d+1})$.
 We have $\alpha-2\beta<0$ which means that $g_{\la, \alpha} (h)\sim \la^{\alpha}u_{\la,\alpha}^{-2}h\to 0$. Consequently, in view of
 \eqref{eq:equivintermpn},  we see that $\varphi^{(\la,\alpha)}(h)\to h^{ (d-1)/2 }$, which implies \eqref{generalbounda}.

For the upper bound of $\varphi^{(\la,\alpha)}(h)$, we use the first term in the minimum in the right-hand side of \eqref{eq:boundspn}. For its lower bound when $h\in [\frac{1}{4}\kappa_d^{-2/ (d+1) },\frac{1}{2}u_{\la,\alpha}^2]$, we notice that $\la^{ \alpha/2}u_{\la,\alpha}^{-1}$ goes to 0, which implies that we need to take the first term of the minimum in the left hand side of \eqref{eq:boundspn}. We then use the fact that the product $h^{(d-1)/2}(1-u_{\la,\alpha}^{-2}h_2)^{3(d-1)/2}$ is bounded from below.

\vskip.2cm

\noindent{$\cdot$ The case  $\alpha \in (\frac{-2}{d-1},0)$.}
The proofs of \eqref{generalbounda} and \eqref{eq:upperboundphi} are identical to the previous case.
For the lower bound $\varphi^{(\la,\alpha)}(h)$, we need again to take the first term of the minimum in the lower bound from \eqref{eq:boundspn}. We obtain
$$
c h^{\frac{d-1}{2}}((1-\frac{1+\la^{-\alpha}}{2}u_{\la,\alpha}^{-2}h_2)\vee 0)^{\frac{3(d-1)}{2}}\le \varphi^{(\la,\alpha)}(h).
$$
To derive \eqref{eq:lowerboundmu}, it then remains to notice that the left-hand side in the above inequality is bounded below by a constant when $h\in [\frac{1}{4}\kappa_d^{-2/(d+1)},\frac{1}{2}\la^{\alpha}u_{\la,\alpha}^2]$.

\vskip.2cm

\noindent{\it (d)  The regime  $\alpha = \frac{-2} {d-1}$.}
For $h\ge 2\kappa_d^{-2/(d+1)}$ we have $h> \un$, which implies that $\stwo(g_{\la, \alpha} (h))=0$ and in turn that $\varphi^{(\la,\alpha)}(h)=0$.  When $h < 2\kappa_d^{ -2/(d + 1)}$, we obtain
$$
g_{\la, \alpha} (h)\sim \la^{-4\beta}\kappa_d^{\frac{2}{d+1}}h(1-\frac{\kappa_d^{\frac{2}{d+1}}h}{2}).
$$
Using  Lemma \ref{lem:proppn}(ii), we obtain that $\varphi^{(\la,\alpha)}(h)$ converges to $\sone(\kappa_d^{ 2/(d+1)}h)$
where $\sone(\cdot)$ is defined at \eqref{eq:defs1}.  Thus \eqref{generalbounda} holds.

The upper bound of $\varphi^{(\la,\alpha)}(h)$ is obtained by taking the first part of the minimum in the upper bound of \eqref{eq:boundspn}. For the lower bound, %We now establish the asserted bounds for $\mu_{n,\alpha}(v,r,h_1,h_2)$.
we need to take the first part of the minimum in the lower bound from \eqref{eq:boundspn}, which gives
$$
c h^{\frac{d-1}{2}}((1-\frac{1+\la^{-\frac{2}{d-1}}}{2}\kappa_d^{\frac{2}{d+1}}h)\vee 0)^{\frac{3(d-1)}{2}}\le \varphi^{(\la,\alpha)}(h).
$$
The left-hand side in the inequality above is bounded from below by a constant as soon as $h\in [\frac{1}{4}\kappa_d^{-2/(d+1)}, \frac{1}{2}\kappa_d^{-2/(d+1)} ].$
% Taking the first part of the minimum in the upper bound of \eqref{eq:boundspn}, we get that $\mu_{n,\alpha}(v,r,h_1,h_2)$ is bounded from above by $c h_2^{\frac{d-1}{2}}$, which is bounded by a constant when $h_2\le \frac{3}{2}\kappa_d^{-2/(d+1)}$. In conclusion, we have that $\mu_{n,\alpha}(v,r,h_1,h_2)$ is bounded from below and from above by positive constants.

\vskip.2cm

\noindent{\it (e)  The regime  $\alpha \in (-\infty,  -\frac{2}{d-1})$.}  We only consider the case $h \in (0,  \un)$ as  otherwise for  $h > \un $ we have $\stwo(g_{\la, \alpha} (h))=0$ and $\varphi^{(\la,\alpha)}(h)=0$. We are in the situation where the support of $\varphi^{(\la,\alpha)}$ shrinks to $\{0\}$ with a bounded total mass.
In view of \eqref{eq:formulapnf}, we claim that in order to show the convergence in distribution \eqref{generalbounda}, it is enough to prove that the integral of $\varphi^{(\la,\alpha)}$ over  $[0,h^{(\la,\alpha)}]$ goes to $1$. To do so, we apply the change of variable
$h = h^{(\la,\alpha)} h'$ to obtain
\begin{equation}
  \label{eq:intalphasmall}
\int_0^{h^{(\la,\alpha)}}\varphi^{(\la,\alpha)}(h)\mathrm{d}h=h^{(\la,\alpha)}\int_0^1\varphi^{(\la,\alpha)}(h^{(\la,\alpha)}h')\mathrm{d}h'.
\end{equation}
We then use the dominated convergence theorem. Indeed, we get from \eqref{eq:defgnh}
$$
g_{\la, \alpha} (h^{(\la,\alpha)}h')\sim 2\la^{2\alpha}h'(1-h').
$$
Inserting this estimate into \eqref{eq:equivintermpn}, we deduce that for all $h'\in (0,1)$ that as $\la \to \infty$
$$
h^{(\la,\alpha)}\varphi^{(\la,\alpha)}(h^{(\la,\alpha)}h')\to \frac{2^d\kappa_{d-1}}{\kappa_d}\left(h'(1-h')\right)^{\frac{d-1}{2}}.
$$
Moreover using the first term in the minimum of the upper bound in \eqref{eq:boundspn} and $h^{(\la,\alpha)}\le \la^{\alpha+\frac{2}{d-1}}$, we obtain
$$h^{(\la,\alpha)}\varphi^{(\la,\alpha)}(h^{(\la,\alpha)}h')\le c.$$
Consequently, thanks to the dominated convergence theorem, we get
\begin{align*}
& \lim_{\la \to \infty} h^{(\la,\alpha)} \int_0^1 \varphi^{(\la,\alpha)}(h^{(\la,\alpha)}h')\mathrm{d}h' \\
& =  \frac{2^{d } \kappa_{d-1}}{\kappa_d}\int_0^1\left(h(1-h')\right)^{\frac{d-1}{2}}\mathrm{d}h'=\frac{2^{d } \kappa_{d-1}}{\kappa_d}B(\frac{d+1}{2},\frac{d+1}{2}),
\end{align*}
where for all $x, y > 0$ we recall that  $ B(x,y):= \int_0^1 t^{x-1} (1 - t)^{y-1} \mathrm{d}t$ is the beta function.

The identity
$$
B(\frac{d+1}{2},\frac{d+1}{2}) = \frac{ \Gamma( \frac{ d+ 1} {2} ) \Gamma( \frac{ d+ 1} {2} ) }  { \Gamma( d + 1) }
$$
and the Legendre duplication formula $\Gamma(2z) =  2^{2z - 1}  \Gamma(z) \Gamma( z + \frac{1} {2} )  / \sqrt{\pi }$
yield
$$
 \frac{2^{ d }\kappa_{d-1}}{\kappa_d}B(\frac{d+1}{2},\frac{d+1}{2}) = 1,
 $$
which shows the required convergence \eqref{generalbounda}.
To obtain \eqref{eq:lowerboundmualphasmall}, it is enough to show that the integral  at \eqref{eq:intalphasmall} is upper and lower bounded by a positive constant. For the upper bound, we integrate the first part of the minimum in the upper bound of \eqref{eq:boundspn} to get
\begin{align*}
& h^{(\la,\alpha)}\int_{0}^{1}
\varphi^{(\la,\alpha)}(h^{(\la,\alpha)}h')\mathrm{d}h'
%\\&
\le c e(\la)\left(h^{(\la,\alpha)}\right)^{\frac{d+1}{2}}\hspace*{-.2cm}(\la^{\frac{\alpha}{2}}u_{\la,\alpha}^{-1})^{d-1}\hspace*{-.2cm}\int_0^1h'^{\frac{d-1}{2}}\mathrm{d}h'\le c.
\end{align*}

For the lower bound, we proceed similarly by integrating the first part of the minimum in the lower bound of \eqref{eq:boundspn}.
This completes the proof of  Lemma \ref{densitybounds}.
%(iv) holds and the proof of Lemma \ref{densitybounds} is complete.
\qed

\vskip.2cm

\noindent{\bf 2.3. Proof of Proposition \ref{prop1}}. The convergence \eqref{Tnconv1} follows directly from  \eqref{generalbounda}.
 The convergence \eqref{Tnconv} follows since the number of points of $T^{(n,\alpha)}(\tilde{\X}_{n,\alpha})$
 belonging to a Borel set $A$ of $\R^{d-1}\times [0, \infty)$ is binomial with mean $\mu^{(n,\alpha)}(A)$.
\qed

\section{ Stabilization and moment bounds}\label{sec:stab}
\allco

Let ${\X}^{(n, \alpha)}  := T^{(n, \alpha)} ( \tilde{\X}_{n, \alpha})$ and ${\P}^{(\la, \alpha)} := T^{(\la, \alpha)}  ( \tilde{\P}_{\la, \alpha})$. Recall we denote the closure of the injectivity region $\{ v \in T(\S^{d-1}), |v| \leq \pi \} $ of the exponential map by $B_{d-1}(\pi)$ where  $T(\S^{d-1})$ is the tangent space at ${\bf u}_0$.

As in \cite{CSY}, for all $k \in \{0, ..., d - 1\}$, finite $\X \subset \R^d$, $x \in \X$, we define the score function $\xi_k(x, \X)$ to be the
product of $(k + 1)^{-1}$ and the  number of $k$-dimensional faces of the convex hull of $\X$ containing $x$; if no faces contain $x$ then
 otherwise we put $\xi_k(x, \X)$ to be zero.

For $x \in \X \subset W_n$ we put
$$
\xi_k^{(n, \alpha)}(x, \X):= \xi_k( [T^{(n, \alpha)} ]^{-1}(x), [T^{(n, \alpha)} ]^{-1}(\X))
$$
and we define $\xi_k^{(\la, \alpha)}(x, \X)$ similarly.   Then
$$
f_k(K_{n, \alpha}) = \sum_{ x \in {\X}^{(n, \alpha)}    }  \xi_k^{(n, \alpha)} ( x, {\X}^{(n, \alpha)}  ) = \sum_{ \tilde{X}_i \in \tilde{\X}_{n, \alpha} }  \xi_k( \tilde{X}_i, \tilde{\X}_{n, \alpha} )
$$
and
\be \label{basicsum}
f_k(K_{\la,\alpha}) = \sum_{x \in  {\P}^{(\la, \alpha)} } \xi_k^{(\la, \alpha)} (x, {\P}^{(\la, \alpha)} ) = \sum_{ x \in  \tilde{\P}_{\la, \alpha} } \xi_k(x,  \tilde{\P}_{\la, \alpha} ).
\ee
Let $u_{\la,\alpha}$ be as in Table \ref{tab:summary} and let $\ul$ be as in \eqref{udef}.   Letting $W_\la$ be the image by $T^{(\la, \alpha)}$ of the region carrying the perturbed points, we have
\be \label{defW}
W_\la :=
u_{\la,\alpha} B_{d-1}(\pi) \times [0, \ul  ]
\ee
We also put
\be \label{defWa}
W_{\infty} :=
\begin{cases}
\R^{d-1} \times \R^+
  &  \alpha \in (\frac{-2} {d - 1}, \infty)\\
\R^{d-1}\times [0,2\kappa_d^{-\frac{2}{d+1}}]& \alpha=-\frac{2}{d-1}\\
\R^{d-1} \times \{\0\}
 & \alpha \in (-\infty,  \frac{-2} { d - 1 }).
\end{cases}
\ee

For $v \in \R^{d-1}$ and $r \in (0, \infty)$ define the cylinder $C_{d-1}(v, r):= B_{d-1}(v, r) \times \R$.  Given a score function $\xi_k^{(\la, \alpha)}, \la \in [1, \infty],$ and $w := (v,h) \in W_\la, \la \in [1, \infty],$ we write
$$
\xi_{k, [r]}^{(\la, \alpha)}(w, {\P}^{(\la, \alpha)}) := \xi_{k}^{(\la, \alpha)}(w, {\P}^{(\la, \alpha)} \cap C_{d-1}(v, r)).
$$
Given $\xi_{k}^{(\la, \alpha)}$,
let $R:= R^{ \xi_{k}^{(\la, \alpha)}}(w, {\P}^{(\la, \alpha)})$ be the smallest integer such that a.s.
$$
%\be \label{defstab}
\xi_{k}^{(\la, \alpha)}(w, {\P}^{(\la, \alpha)}) = \xi_{k, [r]}^{(\la, \alpha)}(w, {\P}^{(\la, \alpha)}) \ \ \forall r \in [R, \infty).
$$
Then $R$ is a random variable \cite{Pe}
 and is a radius of (spatial) stabilization  for $\xi_{k}^{(\la, \alpha)}$ at $w$ with respect to ${\P}^{(\la, \alpha)}$.

In the first part of this section, we show in Lemma \ref{stab} that $R^{ \xi_{k}^{(\la, \alpha)}}(w, {\P}^{(\la, \alpha)})$, $\alpha \in ( - \infty, \infty)$,  have exponentially decaying tails.  When $\alpha \in ( \frac{-2} {d -1}, \infty)$, standard methods show that
 given $w \in W_\la$, we have as $\la \to \infty$
$$
\xi_k^{(\la, \alpha)}(w, {\P}^{(\la, \alpha)} ) \tod \xi_k^{(\infty, \alpha)}(w, \P^{(\infty, \alpha)} ).
$$
When $\alpha \in (-\infty, \frac{-2} {d -1} )$, the limit point process  ${\P}^{(\infty, \alpha)}$ concentrates on a subspace of $\R^{d-1} \times \R$ and, as shown in the proof of Proposition \ref{expectation-1},   we require new  methods  to show the above convergence.
In the second part of this section, we show in Lemma \ref{mombounds} that the rescaled functionals $\xi^{(\la, \alpha)}_k$   satisfy moment bounds of all orders when $\alpha \in  (- \infty, \infty)$. Stabilization and moment bounds yield convergence of expectations.

\begin{prop} \label{expectation}  For all  $\alpha \in [ \frac{-2} {d - 1}, \infty)$ and all $w \in W_\la$ we have
\be \label{LconvPo}
\lim_{\la \to \infty} \E \xi_k^{(\la, \alpha)}(w, {\P}^{(\la, \alpha)} ) =  \E \xi_k^{(\infty, \alpha)}(w, \P^{(\infty, \alpha)}  ).
\ee
We also have for all $w_1, w_2 \in W_\la$
\begin{align} \label{LconvPoisson}
&\lim_{\la \to \infty} \E \xi_k^{(\la, \alpha)}(w_1 , {\P}^{(\la, \alpha)} ) \xi_k^{(\la, \alpha)}(w_2 , {\P}^{(\la, \alpha)} ) \nonumber \\
 & =  \E \xi_k^{(\infty, \alpha)}(w_1, \P^{(\infty, \alpha)}  ) \xi_k^{(\infty, \alpha)}(w_2, \P^{(\infty, \alpha)} ).
\end{align}
\end{prop}

When  $\alpha \in (- \infty, \frac{-2} {d - 1})$ we shall require convergence of one and two point correlation functions holding uniformly over
  $W_\la$.  This result, which relies on coupling arguments, goes beyond standard convergence results for score functions of stabilizing input and is needed to show variance asymptotics.

\begin{prop} \label{expectation-1}
For all  $\alpha \in (- \infty, \frac{-2} {d - 1})$
we have
\be \label{LconvPo-1}
\lim_{\la \to \infty} \E \sup_{(\0,h) \in W_\la} | \xi_k^{( \lambda, \alpha) } ( (\0,h), {\cal P}^{(\la, \alpha)} ) - \xi_k^{( \infty, \alpha) } ( (\0,0), {\cal P}^{(\infty, \alpha)} )| = 0.
\ee
We also have
\begin{align}  \label{LconvPo-3}
& \lim_{\la \to \infty} \E \sup_{(\0,h) \in W_\la} \sup_{(v_1,h_1) \in W_\la} | \xi_k^{( \lambda, \alpha) } ( (\0,h), {\cal P}^{(\la, \alpha)} )
 \xi_k^{( \lambda, \alpha) } ( (v_1,h_1), {\cal P}^{(\la, \alpha)} ) \nonumber  \\
&  \hskip2cm - \xi_k^{( \infty, \alpha) } ( (\0,0), {\cal P}^{(\infty, \alpha)} ) \xi_k^{( \infty, \alpha) } ( (v_1, 0), {\cal P}^{(\infty, \alpha)} ) | = 0.
\end{align}
\end{prop}

The proof of Proposition \ref{expectation} is standard (cf. \cite{BY05, Pe07}) and will be omitted  whereas the proof of Proposition \ref{expectation-1} involves  new ingredients and is deferred to Subsection 3.3.

\vskip.2cm

\noindent{\bf 3.1. Stabilization of score functions.}
We provide rigorous yet straightforward and self-contained arguments showing that $\xi_k^{(\la, \alpha)}, \ k \in \{0,...,d-1 \}$,
$\alpha \in [\frac{-2} {d -1}, \infty)$, have radii of stabilization with exponentially decaying tails.   We shall consider the finite-size scaling images of balls and half-spaces under $T^{(\la, \alpha)}$.
 Given $\X \subset B_d(\0, 1 + \la^{\alpha})$ a finite point set, a point $x_0 \in \X$
 is a vertex of the convex hull of $\X$ iff $B_d(\frac{x_0} {2}, \frac{|x_0|} {2})$ is not covered by $\cup_{x \in \X, x \neq x_0} B_d(\frac{x} {2}, \frac{|x|} {2}).$  The transformation $T^{(\la, \alpha)}$ preserves this geometric characterization of extreme points, which is seen as follows.

 Let $h = u_{\la,\alpha}^2 (1 - \frac{ |x| } { 1 + \la^{\alpha} } )$ and $h_0 := u_{\la,\alpha}^2 (1 - \frac{ |x_0| } { 1 + \la^{\alpha} } )$.  As in \cite{CSY}, the map $T^{(\la, \alpha)} $ at  \eqref{map2} sends  $B_d( \frac{x_0} {2}, \frac{|x_0|} {2})$ to the solid
\begin{align} \label{up1}
& [ \Pi^{\uparrow} (v_0, h_0)]^{(\lambda, \alpha)}  \nonumber \\
& \ \ \ := \{ (v,h) \in W_\la: h \geq u_{\la,\alpha}^2  (1 - \cos (e_\la(v, v_0))) + h_0 \cos (e_\la(v, v_0)) \},
\end{align}
where for all $v_1, v_2 \in u_{\la,\alpha} B_{d-1}(\pi)$ we have
$$
e_\la(v_1, v_2):= d_{ \mathbb{S}^{d-1} }( \exp_{d - 1} (u_{\la,\alpha}^{-1} v_1), \  \exp_{d - 1} (u_{\la,\alpha}^{-1} v_2) )
$$
where $d_{ \mathbb{S}^{d-1} }$ is the geodesic distance on $\mathbb{S}^{d-1}$.
Put
$$
[ \Pi^{\uparrow} (v_0, h_0)]^{(\infty)} := \{ (v,h) \in \R^{d-1} \times \R: h \geq h_0  + \frac{ | v - v_0|^2 } {2}   \}.
$$
Given the cylinder $C_{d-1}(v_0, L), L \in (0, \infty)$, the graph of the boundary of $[ \Pi^{\uparrow} (v_0, h_0)]^{(\la, \alpha)} \cap C_{d-1}(v_0, L)$ converges as $\la \to \infty$ to the graph of the boundary of $[ \Pi^{\uparrow} (v_0, h_0)]^{(\infty, \alpha)} \cap C_{d-1}(v_0, L)$.  Consequently,  $[ \Pi^{\uparrow} (v_0, h_0)]^{(\la, \alpha)}, \la \in [1, \infty),$ are nearly parabolic in a neighborhood of $(v_0, h_0)$ and hence we call them
`quasi up-paraboloids'.  Thus, if $x_0 \in \X$ is sent to $w_0$, we see that $x_0$ is extreme iff $[\Pi^{\uparrow} (w_0)]^{(\la, \alpha)}$ is not covered by $$
\bigcup_{w \in T^{(\la, \alpha)} (\X), \ w \neq w_0 } [\Pi^{\uparrow} (w)]^{(\la, \alpha)}.
$$

The image of the half-space
$$
 \{ x \in \R^d:  \langle x, \frac{x_0} {|x_0|} \rangle  \geq |x_0| \}
$$
 under  $T^{(\la, \alpha)}$ is found by interchanging $h_0$ and $h$ in the right-hand side of \eqref{up1}.  This gives a dual image, namely the `quasi down-paraboloid'
\begin{align} \label{down1}
  [ \Pi^{\downarrow} (v_0, h_0)]^{(\lambda, \alpha)} & := \{ (v,h) \in W_\la: h \cos (e_\la(v, v_0)) \nonumber  \\
 & \ \ \ \ \ \ \ \ \ \ \ \ \ \ \ \   \leq h_0 - u_{\la,\alpha}^2 (1 - \cos (e_\la(v, v_0)))  \}.
\end{align}
The solids  $[ \Pi^{\downarrow} (v_0, h_0)]^{(\la, \alpha)}$ are nearly parabolic in a neighborhood of $(v_0, h_0)$ and their boundaries contain the image of facets of $\partial K_{\la,\alpha}.$ Put
$$
[ \Pi^{\downarrow} (v_0, h_0)]^{(\infty)} := \{ (v,h) \in \R^{d-1} \times \R: h \leq h_0  - \frac{ | v - v_0|^2 } {2}   \}.
$$

\vskip.2cm

 The quasi up- and down-paraboloids have a spatial diameter which is well controlled.

\begin{lemm} \label{growth}  For all $\alpha \in ( -\infty, \infty)$  and all $\la \in [1, \infty]$ we have:

\noindent (i) If $(v_1,h_1) \in \partial ( [ \Pi^{\uparrow} (\0, h_0)]^{(\lambda, \alpha)})$ then $|v_1| \leq 2 \sqrt{2}  \sqrt{h_0 + h_1}.$

\noindent (ii) If $(v, 0) \in \partial ( [ \Pi^{\downarrow} (v_0, h_0)]^{(\lambda, \alpha)})$ then for $|v| > |v_0|$ we have $|v| \leq 2 \sqrt{2} \sqrt{h_0} + |v_0|$.
\end{lemm}

\noindent{\em Proof.} We will only consider the case $\la  \in [1, \infty)$;  the case $\la = \infty$ is easier, since in this instance the geometry involves actual paraboloids and not quasi-paraboloids.

\noindent (i) Recall \eqref{up1}. % we have
Now $e_\la(v, \0) = u_{\la,\alpha}^{-1} |v| \in [0, \pi]$.  Let $(v_1,h_1)$ be on the boundary of $[ \Pi^{\uparrow} (\0, h_0)]^{(\lambda, \alpha)}$.  Then
$$
(h_0 + h_1) \geq u_{\la,\alpha}^2  (1 - \cos (e_\la(v_1, \0))) \geq u_{\la,\alpha}^2  \frac{ (e_\la(v_1, \0))^2} {8} = {|v_1|^2 \over 8},
$$
since $1 - \cos \theta \geq \theta^2/8$ when $|\theta| \leq \pi$.

\vskip.2cm
\noindent (ii)
We set $h = 0$ in \eqref{down1} to find the points $(v, 0)$
where the boundary of $[ \Pi^{\downarrow} (v_0, h_0)]^{(\lambda, \alpha)}$ meets the `lower' boundary of $W_\la$.
This gives $h_0 \geq u_{\la,\alpha}^2  (1 - \cos (e_\la(v, v_0)))$.  The inequality $ 1 - \cos \theta \geq \theta^2/8$, $|\theta| \leq \pi$, gives
\begin{align*}
h_0 & \geq  { u_{\la,\alpha}^2   (e_\la(v, v_0))^2 \over 8}  \\
&  \geq  { u_{\la,\alpha}^2    ( e_\la(v, \0) - e_\la(v_0, \0)  )^2 \over 8}   \\
&  \geq { u_{\la,\alpha}^2   ( ( |v| - |v_0|)/ u_{\la,\alpha}   )^2\over 8}
\end{align*}
by the triangle inequality.
\qed

\vskip.2cm

 \begin{lemm} \label{height}  If $\alpha \in (-\infty,  \infty), \la \in [1, \infty]$  and $w_1 := (v_1,h_1) \in \partial([ \Pi^{\uparrow} (\0, h_0)]^{(\lambda, \alpha)})$ is such that  $[ \Pi^{\downarrow} (w_1)]^{(\lambda, \alpha)}) \cap C_{d-1}(\0,t)^c \neq \emptyset$, then $h_1 \geq \frac{t^2} {64}.$
  \end{lemm}

\noindent{\em Proof.} Since $h_1 \geq h_0$, we may assume $h_0 \leq \frac{t^2} {64}$.
In order that $[ \Pi^{\downarrow} (w_1)]^{(\lambda, \alpha)}) \cap C_{d-1}(\0,t)^c \neq \emptyset$, there must be some point $(v, 0)$ on the boundary of $[\Pi^{\downarrow} (w_1)]^{(\lambda, \alpha)})$ with $|v| \in [t, \infty)$.  By Lemma \ref{growth}(ii) we know that
$$
|v| \leq 2  \sqrt{2 h_1}  + |v_1|.
$$
Also, since $(v_1, h_1)$ belongs to the boundary of $[ \Pi^{\uparrow} (\0, h_0)]^{(\lambda, \alpha)}$, it follows from
Lemma \ref{growth}(i) that
$$
|v_1| \leq 2  \sqrt{2} \sqrt{h_0 + h_1}.
$$
To insure $|v| \in [t, \infty)$ we require  $ 2 \sqrt{2} ( \sqrt{h_1} +  \sqrt{h_0 + h_1})  \in [t, \infty)$.  We must have either
$ \sqrt{2} \sqrt{h_1} \geq \frac{t} {4}$ or $ \sqrt{2}\sqrt{h_0 + h_1}  \geq \frac{t} {4}$.  Either case gives $h_1 \geq \frac{t^2}{64}.$  \qed

\vskip.2cm

\begin{lemm} \label{stab} (i) When $\alpha \in ( \frac{-2} {d - 1},  \infty)$ there is a constant $c \in (0, \infty)$  such that for all $k \in \{0,...,d-1 \}$, $\la \in [1, \infty]$,  and $h_0 \in [0, \ul]$ we have
\be \label{stabbound1}
\PP( R^{\xi_k^{(\la, \alpha)}} ((\0, h_0),  {\P}^{(\la, \alpha)}) > t) \leq c \exp( - \frac{t^{d + 1}} {c} ).
\ee
(ii) When $\alpha \in (- \infty, \frac{-2} {d - 1}]$ there is a constant $c \in (0, \infty)$ such that for all
\vskip.1cm
\noindent $k \in \{0,...,d-1 \}$, $\la \in [1, \infty]$, and $h_0 \in [0,  \ul]$  we have
\be \label{stabbound2}
\PP( R^{\xi_k^{(\la, \alpha)}} ((\0, h_0),   {\P}^{(\la, \alpha)}) > t) \leq c \exp( - \frac{t^{d - 1}} {c} ).
\ee
Moreover, for any $v_1 \in u_{\la,\alpha} B_{d-1}(\pi)$, we have for all $h_1 \in [0, \ul ]$  that
 $$R^{\xi_k^{(\la, \alpha)}} ((v_1, h_1), {\P}^{(\la, \alpha)}) \leq R^{\xi_k^{(\la, \alpha)}} ((v_1, \ul),  {\P}^{(\la, \alpha)}).$$
%\vskip.2cm
%\noindent(iii) Bounds identical to those at \eqref{stabbound2}  hold for $\PP( R^{\xi_k^{(n, \alpha)}} ((\0, h_0),   {\X}^{(n, \alpha)}) > t).$
\end{lemm}

\noindent{\em Proof.}
(i) First  consider the case $\alpha \in (\frac{-2} {d - 1}, \infty)$.
Write $\{ R^{\xi_k^{(\la, \alpha)}} ((\0, h_0)) > t \} \subset E_1 \cup E_2$, where
$$
E_1 := \{ R^{\xi_k^{(\la, \alpha)}} ((\0, h_0)) > t, \  (\0, h_0) \notin {\rm {Ext}} ( {\cal P}^{(\la, \alpha)} )     \}
$$
and
$$
E_2:= \{ R^{\xi_k^{(\la, \alpha)}} ((\0, h_0)) > t, \  (\0, h_0) \in {\rm {Ext}} ( {\cal P}^{(\la, \alpha)} )  \}.
$$
We bound $\PP(E_1)$ and $\PP(E_2)$ separately as follows.  %Without loss of generality we may assume $h$ is larger than some finite scalar, say $2$.

If $E_1$ occurs then there is $
w_1 := (v_1,h_1) \in \partial([ \Pi^{\uparrow} (\0, h_0)]^{(\lambda, \alpha)})
$
which is covered by paraboloids with apices in ${\cal P}^{(\la, \alpha)}$ but not by paraboloids with apices in ${\cal P}^{(\la, \alpha)} \cap C_{d-1}(\0,t)$.  In other words, we have
\be \label{downempty}
[\Pi^{\downarrow}(w_1)]^{(\lambda, \alpha)} \cap C_{d-1}(\0, t) \cap {\cal P}^{(\lambda, \alpha)} = \emptyset
\ee
as well as
\be \label{second}
[\Pi^{\downarrow}(w_1)]^{(\lambda, \alpha)} \cap C_{d-1}(\0, t)^c \neq \emptyset.
\ee

By \eqref{eq:lowerboundmu}, the density of the intensity measure $\mu^{(\la,\alpha)}$  is bounded away from zero on $\frac{1}{2}u_{\la,\alpha} B_{d-1}(\pi) \times [\frac{1}{4}\kappa_d^{-2/(d+1)}, \frac{1}{2}u_{\la,\alpha}^2]$.  Thus $$\mu^{(\la,\alpha)}\left([\Pi^{\downarrow}(w_1)]^{(\lambda, \alpha)} \cap C_{d-1}(\0, t)\right)$$ is bounded below by a constant multiple of $\Vol ([\Pi^{\downarrow}(w_1)]^{(\lambda, \alpha)} \cap C_{d-1}(\0, t))$.

Let us show that this volume is at least $c t^{d+1}$: we start by noticing that in view of \eqref{second}, it follows by Lemma \ref{height} that $h_1 \geq \frac{t^2} {64}.$ We now consider two cases.\\
\noindent{$\cdot$} If the maximal height of the intersection of  $[\Pi^{\downarrow}(w_1)]^{(\lambda, \alpha)}$  and  $\partial\C_{d-1}(\0, t)$ exceeds $\frac{h_1} {2}$, then $[\Pi^{\downarrow}(w_1)]^{(\lambda, \alpha)} \cap C_{d-1}(\0, t)$ contains a cylinder of radius $c t$ and of height $c h_1$, i.e. at least $c t^2$. This gives a volume of at least $c t^{d+1}$.\\
\noindent{$\cdot$} If the maximal height of the intersection of  $[\Pi^{\downarrow}(w_1)]^{(\lambda, \alpha)}$  and  $\partial\C_{d-1}(\0, t)$ is at most $\frac{h_1} {2}$, then the subset of
$[\Pi^{\downarrow}(w_1)]^{(\lambda, \alpha)} \cap C_{d-1}(\0, t)$ having height coordinate at least $\frac{h_1} {4}$
contains a translate of $[\Pi^{\downarrow}(w_1')]^{(\lambda, \alpha)}$ with $w_1'=(v_1, \frac{h_1} {4})$. Consequently, its volume is of order $c h_1^{(d+1)/2}$ so is at least $c t^{d+1}$.

Thus
$$
\PP( [\Pi^{\downarrow}(w_1)]^{(\lambda, \alpha)} \cap C_{d - 1} (\0, t)  \cap {\cal P}^{(\lambda, \alpha)} = \emptyset)
\leq c \exp( - \frac{t^{d + 1}} {c} ).
$$

Discretizing $ \partial([ \Pi^{\uparrow} (\0, h_0)]^{(\lambda, \alpha)}) \cap (\R^{d-1} \times [\frac{t^2} {64}, \infty))$  we obtain
$$
\PP(E_1)  \leq c \exp( - \frac{t^{d + 1}} {c} ).
$$

Now we bound $\PP(E_2)$. If $E_2$ occurs then there is a $w_1:= (v_1, h_1) \in C_{d-1}(\0,t)^c \cap [ \Pi^{\uparrow} (\0, h_0)]^{(\lambda, \alpha)}$ such that $w_1$ is not covered by quasi up-paraboloids. In other words, $[\Pi^{\downarrow}(w_1)]^{(\lambda, \alpha)} \cap {\cal P}^{(\lambda, \alpha)}  $ must be empty.  The  probability that
 $[\Pi^{\downarrow}(w_1)]^{(\lambda, \alpha)}$ does not contain points of ${\cal P}^{(\lambda, \alpha)}$ decays like $\exp(-c h_1^{(d + 1)/2} )$. Discretizing
 $$
 \partial( [ \Pi^{\uparrow} (\0, h_0)]^{(\lambda, \alpha)} ) \cap ( \R^{d-1} \times [\frac{t^2} {64}, \infty))
 $$
 shows that
 $$
 \PP(E_2)  \leq c \exp( - \frac{ t^{d + 1}} {c}  ).
 $$
Combining the bounds for $\PP(E_1)$ and  $\PP(E_2)$  gives  \eqref{stabbound1} as desired.

\vskip.2cm

\noindent (ii) When $\alpha  \in (- \infty, \frac{-2} {d - 1}]$  we follow the same reasoning, i.e. we introduce the events $E_1$ and $E_2$ and estimate $\mu^{(\la,\alpha)}\left([\Pi^{\downarrow}(w_1)]^{(\lambda, \alpha)} \cap C_{d - 1} (\0, t)\right)$ and $\mu^{(\la,\alpha)}\left([\Pi^{\downarrow}(w_1)]^{(\lambda, \alpha)}\right)$. To obtain \eqref{stabbound2}, we only need to show that this measure is at least $c t^{d-1}$ and to do so, we consider two cases.
\\
\noindent{$\cdot$} If the maximal height of the intersection of
$[\Pi^{\downarrow}(w_1)]^{(\lambda, \alpha)}$ and  $\partial\C_{d-1}(\0, t)$ exceeds $\ul$, then $[\Pi^{\downarrow}(w_1)]^{(\lambda, \alpha)} \cap C_{d-1}(\0, t)$ contains a cylinder $B_{d-1}(v_1,ct)\times (0, \ul)$. By \eqref{eq:lowerboundmu} and \eqref{eq:lowerboundmualphasmall}, the measure of such a cylinder is bounded from below by a multiple of the $(d-1)$-dimensional Lebesgue measure of its base, which is of order $c t^{d-1}$.\\

\noindent{$\cdot$} If the maximal height of the intersection of $[\Pi^{\downarrow}(w_1)]^{(\lambda, \alpha)}$ and  $\partial\C_{d-1}(\0, t)$ is at most $\ul$, then we claim that $[\Pi^{\downarrow}(w_1)]^{(\lambda, \alpha)} \cap C_{d-1}(\0, t)$ contains a cylinder $B_{d-1}(v_1',r)\times (0, \ul)$ where $v_1'$ belongs to $B_{d-1}(v_1,t)$ and $r$ is proportional to $\sqrt{h_1}$. Indeed, if $v$ is such that $(v, \ul)$ belongs to $\partial\,[\Pi^{\downarrow}(w_1')]^{(\lambda, \alpha)}$, then we get from \eqref{down1} that
$$
\cos(e_\lambda(v,v_1))=1-\frac{h_1- \ul} {u_{\la,\alpha}^2 - \ul}.
$$
In particular, this implies that for large $\la$, $|v-v_1|\sim \sqrt{2h_1}$ and in turn, we obtain that $[\Pi^{\downarrow}(w_1)]^{(\lambda, \alpha)} \cap C_{d-1}(\0, t)$ contains a cylinder of radius $c \sqrt{h_1}$, which is at least $c t$ in view of \eqref{second} and Lemma \ref{height}. Again by Lemma \ref{densitybounds}(iv), the measure of such a cylinder is bounded below by $c t^{d-1}$.
This completes the proof of Lemma \ref{stab}.  \qed

%\noindent (iii) Straightforward modifications of the above proofs show that \eqref{stabbound1} and \eqref{stabbound2} hold when $\P^{(\la, \alpha)}$ is %replaced by $\X^{(n, \alpha)}$.
%This completes the proof of Lemma \ref{stab}.  \qed

\vskip.5cm

\noindent {\bf 3.2. Moment bounds for score functions.}  We next establish that moment bounds for
 score functions hold uniformly in the spatial coordinates.
\vskip.2cm

\begin{lemm} \label{mombounds} (i) For all $\alpha \in (\frac{-2} {d - 1}, \infty)$,  $p \in [1, \infty)$, $k \in \{0,1,..., d-1 \}$ and $d = 1,2,...$ there is a constant $c:= c(\alpha, p,k,d) \in (0, \infty)$  such that
$$
 \sup_{1 \leq \la \leq \infty} \sup_{v \in u_{\la,\alpha} B_{d-1}(\pi)}  \E ( \xi_k^{(\la, \alpha)} ((v,h), {\P}^{(\la, \alpha)} ) )^p  \leq c h^{c } \exp( -  \frac{ h^{(d + 1)/2}} {c} ).  \label{Poissbound}
$$
(ii) If  $\alpha  \in (- \infty,  \frac{-2} {d - 1}]$ then similar bounds hold, with $\exp( -  \frac{ h^{(d + 1)/2}} {c} )$ replaced by
$\exp( -  \frac{ h^{(d -1)/2}} {c} )$. %and $h^{(\red{d + 3}) pk/2 }$ replaced by $c^{pk}$, $c$ a constant.
\end{lemm}

\vskip.2cm

\noindent{\em Proof.}  (i)  Let $w:=(v,h) \in W_\la$ be fixed. Let $H:= H(w)$ be the maximal height coordinate with respect to $\R^{d - 1}$ of all apices of quasi-down paraboloids $[ \Pi^{\downarrow} (v_1, h_1)]^{(\la, \alpha)}$ containing a  face of the boundary of $T^{(\la, \alpha)}( K_{\la,\alpha})$ and also containing $w$; if no such quasi-down paraboloid exists then we put $H(w)$ to be zero. By \eqref{eq:lowerboundmu}, for all  $\alpha \in ( -\frac{2} {d - 1}, \infty)$,
the density of the intensity measure of $\P^{(\la, \alpha)}$ is bounded away from zero for heights $h \in [1, \frac{u_{\la,\alpha}^2} {2}]$. Since the volume of $[ \Pi^{\downarrow} (v_1, h_1)]^{(\la, \alpha)}$ is proportional to $h_1^{(d + 1)/2}$ it follows that for all $w \in W_\la$
\be \label{Height}
\PP(H(w) \geq t) \leq c \exp( - \frac{t^{(d + 1)/2 } } {c} ),  \ t > 0.
\ee

Abbreviate the radius of stabilization $R^{ \xi_k^{(\la, \alpha) } }( w, \P^{(\la, \alpha)})$ by $R$. Let $N(w)$ be the
cardinality  of the point set $ {\P}^{(\la, \alpha)} \cap (B_{d-1}(v, R) \times [0, H]).$
For all $k \in \{0,...,d-1\}$ we have
\be \label{generalbound}
\xi_k^{(\la, \alpha)} (w, {\P}^{(\la, \alpha)} ) \leq \frac{1} { k + 1}  \binom{N(w)} {k - 1} {\bf {1} } ( w \ \rm{is} \ \rm{extreme} ).
\ee

By \eqref{eq:upperboundphi}, the density of $\mu^{(\la, \alpha)}$ is upper bounded by $c_2h^{(d-1)/2}$ on cylinders of height $h$, giving
\be \label{contentbound}
\mu^{(\la, \alpha)}  ( B_{d - 1}(v, r) \times [0, h] ) \leq c \Vol_{d-1}( B_{d-1}(v, r)) \times h^{(d + 1)/2},
\ee
where $\Vol_{d-1}(A)$ denotes the Lebesgue measure of a full $(d-1)$-dimensional subset $A$ of $\R^{d-1}$.
We obtain as in the proof of Lemma 4.4 of \cite{CY2}
$$
\E N^{pk} \leq c \cdot (h + 1)^{c},
$$
where $c:= c(\alpha, p,k,d) \in (0, \infty)$.

The above bounds (with $t$ set to $h$ in \eqref{Height}) together with the Cauchy-Schwarz inequality
$$
\E ( \xi_k^{(\la, \alpha)} (w, {\P}^{(\la, \alpha)} )) ^p \leq (\E (\xi_k^{(\la, \alpha)} (w, {\P}^{(\la, \alpha)} ))^{2p})^{1/2} \PP(  \xi_k^{(\la, \alpha)} (w, {\P}^{(\la, \alpha)} ) > 0)^{1/2}
$$
yield the claimed bound.

 \vskip.2cm

 \noindent (ii) When $\alpha \in (-\infty, \frac{-2} {d -1} ]$, we follow the proof of part (i) nearly verbatim, save for the following modifications.
  %Use \eqref{bounds-6} and
  Use the method of proof of Lemma \ref{stab}(ii)  to show for all $w \in W_\la$ that
\be \label{Height-1}
\PP(H(w) \geq t) \leq c \exp( - \frac{t^{(d - 1)/2 } } {c} ),  \ t > 0.
\ee
By Lemma \ref{densitybounds}(iii),  the bound \eqref{contentbound} is replaced by
 \be \label{contentbound-1}
\mu^{(\la, \alpha)}  ( B_{d - 1}(v, r) \times [0, h] ) %\leq c \cdot  \mathrm{d} \H^{(\la, \alpha)} (B_{d-1}(v, r))
\leq c \cdot \Vol_{d-1} (B_{d-1}(v, r)),
\ee
which holds for all $h \in [0, \ul ]$.  With these small modifications, we may now follow the proof of part (i).
   \qed

\vskip.2cm

\noindent{\bf 3.3. Proof of Proposition \ref{expectation-1}.}

\noindent{\em (i) Proof of \eqref{LconvPo-1}.}
 Given finite point sets ${\cal X}$ and ${\cal \Y}$ in $\R^{d-1} \times \R^+$, we let ${\bf d} ({\cal X}, {\cal Y} ) := \sup_{x \in \cal \X} \inf \{ |x - y|: \ y \in \cal \Y)$.
%For $r, l > 0$ and $v \in \R^{d-1} \times \{0 \}$ we put
%\begin{equation} \label{defSvrl}
%S(v, r, l):= B_{d-1}(v, r) \times [0, l].
%\end{equation}

We claim that the scores $\xi_k^{( \lambda, \alpha) }$,  $k \in \{1,...,d \}$, are continuous in the sense that if ${\cal X} \subset \R^{d-1} \times \{0 \}$ is in regular position (that is to say that for all $k \in \{1,...,d \}$ the intersection of the boundaries of $k$ quasi paraboloids contains at most $(d - k + 1)$ points of $\X$),  if
 ${\bf d}({\cal X},  {\cal X}_\la) \to 0$ as $\la \to \infty$, where $\X_\la \subset \R^{d-1} \times \R^+$ are finite for all $\la$,
 and if $x_\la \to x$ as $\la \to \infty$ with $x \in \R^{d-1} \times \{0 \}$, then for all $r > 0$ we have
\be \label{cont}
\xi_k^{( \lambda, \alpha) }(x_\la, {\cal X}_\la \cap C_{d-1}(x,r))  \to \xi_k^{( \infty, \alpha) } (x, {\cal X} \cap  C_{d-1}(x,r) ).
\ee
We may see this as follows.  Let $\epsilon > 0$ be the minimal distance between any down-paraboloid containing $d$ points of $\X$ and the rest of the point set in $C_{d-1}(x,r)$.  Perturbations of the paraboloids within an $\epsilon$ parallel set do not change the number of $k$-dimensional faces.  In particular, for $\la$ large enough, the boundary of $\cup_{w \in \X} [ \Pi^{\downarrow} (w) ]^{(\la)} $ is included in that parallel set.  Hence the number of $k$-dimensional faces does not change.  Thus for $\la$ large enough we have for all $k \in \{1,...,d \}$
$$
\xi_k^{( \lambda, \alpha) }(x, {\cal X} \cap C_{d-1}(x,r) ) = \xi_k^{( \infty, \alpha) } (x, {\cal X} \cap C_{d-1}(x,r) ).
$$
Also, since $\cal X$ is in regular position, perturbations of points in $\cal X$ do not change the value of $\xi_k^{( \lambda, \alpha) }(x, {\cal X} \cap C_{d-1}(x,r))$.
Thus if ${\bf d} ({\cal X}_\la, {\cal X}) \to 0$ and $x_\la \to x$ as $\la \to \infty$,  then for $\la \geq \la_0$ large enough  we have
$$
\xi_k^{( \lambda, \alpha) }(x_\la, {\cal X}_\la  \cap C_{d-1}(x,r) ) = \xi_k^{( \infty, \alpha) } (x, {\cal X} \cap C_{d-1}(x,r) ).
$$
This gives the assertion \eqref{cont}.

 Notice that if $y_\la, \la \geq 1,$ is a sequence of points such that $|y_\la - x| \leq |x_\la - x|$ then we also have
for $\la \geq \la_0$
$$
\xi_k^{( \lambda, \alpha) }(y_\la, {\cal X}_\la  \cap C_{d-1}(x,r) ) = \xi_k^{( \infty, \alpha) } (x, {\cal X} \cap C_{d-1}(x,r) ).
$$
This convergence is uniform over all such sequences $y_\la, \la \geq 1,$ that is to say
\be \label{uniform}
\sup_{y_\la:  |y_\la - x| \leq |x_\la - x|}  | \xi_k^{( \lambda, \alpha) }(y_\la, {\cal X}_\la  \cap C_{d-1}(x,r)) - \xi_k^{( \infty, \alpha) } (x, {\cal X} \cap C_{d-1}(x,r))| = 0.
\ee

Next we let ${\cal P}^{(\infty, \alpha)}$ be a Poisson point process on $\R^{d-1} \times \{\0 \}$ of intensity measure with density $1.$  Rescale its restriction to $u_{\la,\alpha} B_{d-1}(\pi)$ by $u_{\la,\alpha}^{-1}$ to get a Poisson point process on
$B_{d-1}(\pi)$ of intensity measure with density $ \frac{\la}{d \kappa_d}.$  Using the exponential map $\exp_{d-1}$ we generate a Poisson point process on $\S^{d-1}$ of intensity measure with density $\frac{\la}{d \kappa_d}.$  Without loss of generality, we may let $\P_\la$ be this Poisson point process and we let $\tilde{\P}_{\la, \alpha}$ be the perturbed point process arising from  $\P_\la$.  %Given $\tilde{\P}_{\la_0, \alpha}$ we generate $\tilde{\P}_{\la, \alpha}$ for $\la \geq \la_0$ by keeping the perturbations $e_i$ associated with $\tilde{\P}_{\la_0, \alpha}$ and adjoining perturbations for the points in $\tilde{\P}_{\la, \alpha} \setminus \tilde{\P}_{\la_0, \alpha}$.

By construction we get $T^{(\la, \alpha)}(\P_\la)$ a.s. converges  to ${\cal P}^{(\infty, \alpha)}$.
Notice that ${\bf d} ( {\cal P}^{(\la, \alpha)},  {\cal P}^{(\infty, \alpha)})$ is a.s. of the order $O( \la^{\alpha}  u_{\la,\alpha}^{2} )= o(1)$ since the spatial coordinates of the points of ${\cal P}^{(\la, \alpha)}$
 are within $\la^{\alpha} u_{\la,\alpha} = o(1)$ of the spatial coordinates of the points of ${\cal P}^{(\infty, \alpha)}$ , whereas the height coordinates of the points of ${\cal P}^{(\la, \alpha)}$ are within $\la^{\alpha} u_{\la,\alpha}^{2}= o(1)$ of the height coordinates of the points of
 ${\cal P}^{(\infty, \alpha)}$.  Thus ${\cal P}^{(\la, \alpha)}$  a.s. converges  to ${\cal P}^{(\infty, \alpha)}$.

In view of \eqref{uniform}, this all implies that, given ${\cal P}^{(\la, \alpha)}(\omega), \la \geq 1,$  if $\la \geq \la_0(\omega)$ is large enough then for a fixed $v \in \la^{\beta} B_{d-1}(\pi)$
\begin{align*}
& \sup_{h: \ (v,h) \in W_\la} | \xi_k^{( \lambda, \alpha) } ( (v,h), {\cal P}^{(\la, \alpha)}  \cap C_{d-1}(v,r)) \\
  & \ \ \ \ \ \ \ \ \ \ \ \ \ \ - \xi_k^{( \infty, \alpha) } ( (v,0), {\cal P}^{(\infty, \alpha)}  \cap C_{d-1}(v,r))| = 0.
\end{align*}
Thus putting $v = \0$ we a.s. obtain
\begin{align}   \label{convergence}
& \lim_{\la \to \infty} \sup_{h: \ (\0,h) \in W_\la} | \xi_k^{( \lambda, \alpha) } ( (\0,h), {\cal P}^{(\la, \alpha)}  \cap C_{d-1}(\0,r)) \nonumber \\
 & \ \ \ \ \ \ \ \ \ \ \ \ \ \ - \xi_k^{( \infty, \alpha) } ( (\0,0), {\cal P}^{(\infty, \alpha)}  \cap C_{d-1}(\0,r))| = 0.
\end{align}
Now for all $p \in [1, \infty)$ the random variables
$$
\sup_{(\0,h) \in W_\la}  \xi_k^{( \lambda, \alpha) } ( (\0,h), {\cal P}^{(\la, \alpha)}  \cap C_{d-1}(\0,r))
$$
have a $p$th moment which is bounded by a constant not depending on $\la$, since the scores $\xi_k^{( \lambda, \alpha) }$
are bounded uniformly in $(\0,h)$ by $\binom{N} {k - 1}$, where $N$ is the cardinality of ${\cal P}^{(\la, \alpha)}  \cap C_{d-1}(\0,r)$ and where  $N$ has moments of all orders.  Likewise the $p$th moment of
$$
  \xi_k^{( \infty, \alpha) } ( (\0,0), {\cal P}^{(\infty, \alpha)}  \cap C_{d-1}(\0,r) )
$$
is bounded by a constant.
Thus we may upgrade \eqref{convergence} to convergence of expectations, i.e.,
\begin{align}  \label{exp-convergence}
& \lim_{\la \to \infty} \E \sup_{(\0,h) \in W_\la} | \xi_k^{( \lambda, \alpha) } ( (\0,h), {\cal P}^{(\la, \alpha)}  \cap C_{d-1}(\0,r) ) \nonumber \\
 & \ \ \ \ \ \ \ \ \ \ \ \ \ - \xi_k^{( \infty, \alpha) } ( (\0,0), {\cal P}^{(\infty, \alpha)}  \cap C_{d-1}(\0,r) )| = 0.
\end{align}

Next, in the regime $\alpha \in (-\infty,  \frac{-2} {d-1} )$ recall that by Lemma \ref{stab}(ii), the radius of stabilization
$$
R:= R^{ \xi_{k}^{(\la, \alpha)}}( (\0, \ul ),  {\P}^{(\la, \alpha)})
$$
is  such that {\em for all} $(\0,h) \in W_\la$ the scores  $\xi_k^{( \lambda, \alpha) } ( (\0,h), {\cal P}^{(\la, \alpha)})$ are determined by data in the cylinder $C_{d-1}(\0,  R).$
In other words, for all $w:= (\0,h) \in W_\la$
$$
\xi_{k}^{(\la, \alpha)}(w, {\P}^{(\la, \alpha)}) = \xi_{k, [r]}^{(\la, \alpha)}(w, {\P}^{(\la, \alpha)}) \ \ \forall r \in [R, \infty).
$$
Since $R$ has exponentially decaying tails standard arguments shows that we may replace the point sets ${\cal P}^{(\la, \alpha)}  \cap C_{d-1}(\0, r)$ and
 ${\cal P}^{(\infty, \alpha)}  \cap C_{d-1}(\0, r)$
  with ${\cal P}^{(\la, \alpha)}$ and ${\cal P}^{(\infty, \alpha)}$, respectively.  We thus
 deduce \eqref{LconvPo-1} from \eqref{exp-convergence}.

\vskip.2cm

\noindent{\em (ii) Proof of \eqref{LconvPo-3}.}  We prove \eqref{LconvPo-3}  by following the methods above.  For all $\la \in [1, \infty]$, $\alpha \in (-\infty,  \frac{-2} {d-1} )$ and $h, h_1, r \in [0, \infty)$,  put
\begin{align*}
& Y_k^{(\la, \alpha)}(0,h,v_1,h_1,r)\\
& := \xi_k^{( \lambda, \alpha) } ( (\0,h), {\cal P}^{(\la, \alpha)}  \cap C_{d-1}(\0,r) )  \xi_k^{( \lambda, \alpha) } ( (v_1,h_1), {\cal P}^{(\la, \alpha)}  \cap C_{d-1}(v,r) ).
\end{align*}
We may  extend \eqref{convergence} to a.s. obtain
\begin{align}  \label{convergence-1}
%& \lim_{\la \to \infty} \sup_{h, h_1: \ (\0,h) \in W_\la, (v_1,h_1) \in W_\la  }
&  | Y_k^{(\la, \alpha)}(0,h,v_1,h_1,r) - Y_k^{(\infty, \alpha)}(0, 0, v_1, 0,r)|  \to 0
\end{align}
uniformly over all $h$ and $h_1$ such that $(\0,h),  (v_1,h_1) \in W_\la$.

Appealing to moment bounds, we may  upgrade \eqref{convergence-1} to get convergence of expectations, i.e.
\begin{align}  \label{exp-convergence-1}
& \lim_{\la \to \infty} \E \sup_{h, h_1: \ (\0,h), (v_1,h_1) \in W_\la  }
 | Y_k^{(\la, \alpha)}(0,h,v_1,h_1,r) - Y_k^{(\infty, \alpha)}(0, 0, v_1, 0,r)| = 0.
\end{align}

Also, in the regime $\alpha \in (-\infty,  \frac{-2} {d-1} )$, given $v_1 \in u_{\la,\alpha}  B_{d-1}(\pi)$,
the scores  $\xi_k^{( \lambda, \alpha) } ( (v_1,h_1), {\cal P}^{(\la, \alpha)})$ are determined by data in $C_{d-1}(v_1, R)$, where
 $$
 R:= R^{ \xi_{k}^{(\la, \alpha)}}( (v_1, \ul ),  {\P}^{(\la, \alpha)}).
 $$
In other words, for all $w_1:= (v_1,h_1) \in W_\la$ we have
$$
\xi_{k}^{(\la, \alpha)}(w_1, {\P}^{(\la, \alpha)}) = \xi_{k, [r]}^{(\la, \alpha)}(w_1, {\P}^{(\la, \alpha)}) \ \ \forall r \in [R, \infty).
$$

Since $R$ has exponentially decaying tails, we may deduce \eqref{LconvPo-3} from \eqref{LconvPo-1}
and \eqref{exp-convergence-1}  via standard arguments.  \qed

\section{ Proofs of main results}\label{sec:proofs}

\allco

The proofs of Theorems \ref{expectasy} and \ref{variance}  require techniques going beyond those in \cite{CSY} and \cite{CY2}.  In part, this is because for $\alpha \in [\frac{-2} {d-1}, \infty)$, the density of the measure $\mu^{(\la, \alpha)}$ explodes for large values of $h$.  Also, for
   $\alpha \in (-\infty, \frac{-2} {d-1} )$   the point processes ${\P}^{(\la, \alpha)}, \la \geq 1,$ and
${\P}^{(\infty, \alpha)}$ concentrate on different spaces, namely $\R^{d-1} \times \R$ and  $\R^{d-1} \times \0$, respectively. We provide a complete proof of Theorem \ref{Thm2}.  Moreover, the proof methods may be easily adapted to establish the proof of Theorem 1.1 of \cite{CY2}, whose proof details were omitted.
  The a.s. convergence of Theorem \ref{Thm3} follows from the coupling introduced in the proof of Proposition \ref{expectation-1}, the techniques used to prove Proposition 5.1 in \cite{CY2}, and the fact that all quasi down-paraboloids are locally nearly parabolic.    Since there are no new essential ideas, we shall leave the details of the proof of  Theorem  \ref{Thm3} to the reader.

\vskip.4cm

\noindent{\bf 4.1. Proof of Theorem \ref{expectasy}. }
Using Proposition \ref{expectation}, we first show expectation asymptotics when binomial input is replaced by Poisson input.  Then we shall de-Poissonize to obtain \eqref{limex2}.
\vskip.1cm
\noindent{\em  Poisson input}.
%We first consider the case $\alpha \in [ \frac{-2} {d-1}, \infty)$.
Recall $\beta:= \beta(\alpha)$, the definition of $W_\la$ at \eqref{defW}, and the definition of  $f_k(K_{\la,\alpha})$ at \eqref{basicsum}.
By the Slivnyak-Mecke formula, with $(v,h)$ standing for a generic point in $W_\la$ we have
\begin{align*}
\E f_k(K_{\la,\alpha}) & =  \int_{B_{d}({\bf 0},1+n^{\alpha})}\E \xi_k(x,\tilde{\P}_{\lambda,\alpha})\mathrm{d}\tilde{\P}_{\lambda,\alpha}(x).
%\int_{u_{\la,\alpha} B_{d-1}(\pi) \times I_{\la, \alpha} } \E \xi_k^{(\la, \alpha)} ((v,h), {\P}^{(\la, \alpha)} ) \mathrm{d} {\P}^{(\la, \alpha)} ((v,h)) \nonumber \\
%& = \Vol( u_{\la,\alpha} B_{d-1}(\pi)  ) \int_{I_{\la, \alpha}}  \E \xi_k^{(\la, \alpha)}((\0,h), {\P}^{(\la, \alpha)} )  \mathrm{d} {\P}^{(\la, \alpha)}_{\0} (h), \label{expect}
\end{align*}
We now apply the change of variables ${\bf u}=\frac{x}{|x|}$ and $h=u_{\la,\alpha}^2(1-\frac{|x|}{1+\lambda^{\alpha}})$. Using Lemma \ref{lem:proppn} (i), we notice that
\begin{equation}
  \label{eq:devPlambdaalpha}
\mathrm{d}\tilde{\P}_{\lambda,\alpha}(x)=u_{\la,\alpha}^{d-1}\mathrm{d}\sigma_{d-1}({\bf u})
\varphi^{(\lambda,\alpha)}(h)\mathrm{d}h.
\end{equation}
Consequently, we get
\begin{align}\label{eq:rewritingEfk}
\E f_k(K_{\la,\alpha}) & =  d\kappa_du_{\la,\alpha}^{d-1}\int_0^{h^{(\la,\alpha)}}\E \xi_k^{(\la, \alpha)}((\0,h), {\P}^{(\la, \alpha)} )
\varphi^{(\lambda,\alpha)}(h)\mathrm{d}h
\end{align}
where we use the rotational invariance of  $\xi_k$, i.e. the fact that $\xi_k(x,\tilde{\P}_{\lambda,\alpha} )$ does not depend on ${\bf u}$ and is equal to $\xi_k^{(\la, \alpha)}((\0,h), {\P}^{(\la, \alpha)} )$.

\vskip.2cm

\noindent (i) {\em The case $\alpha \in [ \frac{-2} {d-1}, \infty)$}.  We subdivide the integration domain in \eqref{eq:rewritingEfk} into
$h\ge \log \la$ and $h\le \log\la$.
The moment bounds of Lemma \ref{mombounds} imply that
 $$
 \E  \xi_k^{(\la, \alpha)}((\0,h), {\P}^{(\la, \alpha)} )  {\bf 1} (h \geq  \log \la )
 $$
 decays faster than $\lambda^{-1}$.
Consequently, in view of  \eqref{eq:totalmassmulambdaalpha} we obtain
\begin{align}\label{eq:negligiblepartEfk}
& d\kappa_du_{\la,\alpha}^{d-1}\int_{\log\la}^{h^{(\la,\alpha)}} \E  \xi_k^{(\la, \alpha)}((\0,h), {\P}^{(\la, \alpha)} )\varphi^{(\lambda,\alpha)}(h)\mathrm{d}h\nonumber \\
&\le o(\la^{-1})u_{\la,\alpha}^{d-1}\int_{0}^{h^{(\la,\alpha)}}\varphi^{(\lambda,\alpha)}(h)\mathrm{d}h=o(1).
\end{align}
We now use the dominated convergence theorem to show that
\begin{align}
  \label{eq:mainpartEfk}
& \lim_{\la\to\infty}\int_0^{\log\la}\E \xi_k^{(\la, \alpha)}((\0,h), {\P}^{(\la, \alpha)} )\varphi^{(\lambda,\alpha)}(h)\mathrm{d}h \nonumber \\
& = \int_0^{\infty}\E \xi_k^{(\infty, \alpha)}((\0,h), {\P}^{(\la, \alpha)} )\mathrm{d}\nu^{(\infty,\alpha)}(h).
\end{align}
By \eqref{eq:convphi} and Proposition \ref{expectation}, the integrand converges a.e. to the desired limit.
By  Lemma \ref{mombounds} and \eqref{eq:upperboundphi},  the integrand is dominated by an exponentially decaying function of $h$ and consequently \eqref{eq:mainpartEfk} holds. Combining it with \eqref{eq:negligiblepartEfk}, we obtain
\be  \label{Poissonlimit1}
\lim_{\la \to \infty} \frac{ \E f_k(K_{\la,\alpha})} { d \kappa_d u_{\la,\alpha}^{d-1} }  =
 \int_{0}^{\infty}  \E \xi_k^{(\infty, \alpha)} ((\0,h), {\P}^{(\infty, \alpha)} )  \mathrm{d} \nu^{(\infty,\alpha)}(h).
\ee

\noindent (ii) {\em The case $\alpha \in (-\infty, \frac{-2} {d-1})$}. This situation is more delicate since the limit point process ${\P}^{(\infty, \alpha)} $  concentrates on a lower-dimensional subset of $\R^{d-1} \times \R^+$ rather than $\R^{d-1} \times \R^+$ itself,
making it difficult to apply standard convergence theorems.
Using \eqref{eq:totalmassmulambdaalpha}, we obtain
\begin{align*}
&\frac{\E f_k(K_{\la,\alpha})}{\la}- \E \xi_k^{(\infty, \alpha)}((\0,0), {\P}^{(\la, \alpha)} )\\& =  \int_0^{h^{(\la,\alpha)}}(\E \xi_k^{(\la, \alpha)}((\0,h), {\P}^{(\la, \alpha)} )-\E \xi_k^{(\infty, \alpha)}((\0,0), {\P}^{(\la, \alpha)} ))\varphi^{(\lambda,\alpha)}(h)\mathrm{d}h.
\end{align*}
Now \eqref{LconvPo-1} implies  that the right-hand side of the above equation goes to zero,
that is to say
\be  \label{Poissonlimit2}
\lim_{\la \to \infty} \frac{ \E f_k(K_{\la,\alpha})} {\la } =  \int_{0}^{\infty}  \E \xi_k^{(\infty, \alpha)} ((\0,h), {\P}^{(\infty, \alpha)} ) \mathrm{d}\nu^{(\infty,\alpha)}(h).
\ee

\vskip.2cm
\noindent{\em  Binomial input}. We use a coupling to de-Poissonize the limits  \eqref{Poissonlimit1} and \eqref{Poissonlimit2}.

\noindent (i) {\em The case $\alpha \in (\frac{-2} {d-1}, \infty)$}.
Enumerate the points in $ {\cal P}^{(n, \alpha)}$ by $Z_1,...,Z_{N(n)}$, $N(n)$ an independent Poisson random variable with parameter $n$ and where the $Z_i$ have distribution $\mu^{(n,\alpha)}$ on $W_n$. % , where we recall the definition of $\mu_{n,\alpha}$ at \eqref{intensity}.
 Consider the coupled point set ${\cal Y}_n$ obtained by discarding or adding i.i.d. points $Z_i$ to $W_n$ where the $Z_i$ have the distribution $\mu^{(n,\alpha)}$:
 \be
 {\cal Y}_n : = \begin{cases} Z_1,...,Z_{N( n)- (N(n) - n )^+}
 , & {\rm if} \ N(n) \geq n \\
Z_1,...,Z_{ N(n) + (n  - N(n))^+}
 , & {\rm if} \ N(n) < n.
\end{cases}
 \ee
Then ${\cal Y}_n \eqd  \{ Z_1, Z_2,..., Z_{n} \} \eqd \X^{(n, \alpha)}$.

Let $\hat{\P}^{( n, \alpha)} := \{ Z_1,...,Z_{N(b_n)} \}$ be the restriction of $ {\cal P}^{(n, \alpha)}$ to $u_{n,\alpha} B_{d-1}(\pi) \times [0, \gamma \log n]$, with $\gamma > 0$ a large constant to be chosen, and where $N(b_n)$ is an independent Poisson random variable with parameter $b_n$, which represents the expected number of points in $u_{n,\alpha}  B_{d-1}(\pi) \times [0, \gamma \log n]$.   By \eqref{generalbound} we have
\be \label{boundfora}
b_n \leq c u_{n,\alpha}^{d-1} \times (\gamma \log n) \times (\gamma \log n)^{ \frac{d + 1} {2}}.
\ee
Let $\hat{\X}^{( n, \alpha)} := \{ Z_1,...,Z_{Bi(n, b_n/n)}\}$ be the restriction of $\X^{(n, \alpha)}$ to $u_{n,\alpha} B_{d-1}(\pi) \times [0, \gamma \log n]$.

For any $\X \subset n^{\beta} B_{d-1}(\pi) \times [0, \gamma \log n]$, define
\be \label{basicsum-1}
F(\X) : = F^{(n, \alpha)}(\X) := \sum_{x \in  \X } \xi_k^{(n, \alpha)} (x, \X ).
\ee
Straightforward  modifications of the above methods show that
%\be  \label{Poissonlimit1-new}
$$
\lim_{\la \to \infty} \frac{ \E F(  \hat{\P}^{( n, \alpha)} ) } { n }  =
  \int_{0}^{\infty}  \E \xi_k^{(\infty, \alpha)} ((\0,h), {\P}^{(\infty, \alpha)} ) \mathrm{d}\nu^{(\infty,\alpha)}(h)
$$
and it suffices to show the same limit when $\hat{\P}^{( n, \alpha)} $ is replaced by $\hat{\X}^{( n, \alpha)}$.
To do this it suffices to  show
\be \label{smalldifference}
\E |  F(  \hat{\P}^{( n, \alpha)} ) -  F(  \hat{\X}^{( n, \alpha)}   )  | = o( u_{n,\alpha}^{d-1} ).
\ee

Put  $E_{n, 1} := \{ N(b_n) \leq {\rm{Bi}}(n, b_n/n)\}$.
We first consider the expected add-one cost
$$
| \E (F(  \hat{\P}^{( n, \alpha)} ) -  F(  \hat{\P}^{( n, \alpha)} \cup Z_{N(b_n) + 1} )) |.
$$
Let $E_{n,2} $ be the event that the radii of stabilization of $\xi^{(n, \alpha)}$ at all points in  $\hat{\P}^{( n, \alpha)}$ are bounded by $\gamma \log n$. If $\gamma$ is large enough then $\PP(E_{n,2} ^c) = o (n^{-d - 1})$.
Conditional on $Z_{N(b_n) + 1} = (v_0, h_0)$, on the event $E_{n,2}$, the number of extreme points in $\hat{\P}^{( n, \alpha)}$ whose score may be modified by the insertion of $Z_{N(b_n) + 1}$ is bounded by the  number of points of $\hat{\P}^{( n, \alpha)}$
in $B_{d-1}(v_0, \gamma \log n) \times [0, \gamma \log n]$.
 Recalling \eqref{generalbound}, it may be shown that there is a third event $E_{n,3}$, $\PP(E_{n,3}^c) = o (n^{-d - 1})$,  such that on $ E_{n,1} \cap E_{n,2} \cap E_{n,3}$  we have
 $$
 {\rm{card} } ( \hat{\P}^{( n, \alpha)} \cap (B_{d - 1} (v_0, \gamma \log n) \times [0, \gamma \log n]) )
   \leq c(\gamma)(\log n)^{d} \times (\log n)^{\frac{d + 1} {2} },
 $$
 where the first factor is the volume of the set carrying the points and the second factor is an upper bound on the intensity measure.

Generously bounding the $f_k$ statistic on an $n$-point set by $cn^d$, we deduce
$$
 \E | ( F(  \hat{\P}^{( n, \alpha)} ) -  F(  \hat{\P}^{( n, \alpha)} \cup Z_{N(b_n) + 1} )) {\bf 1} (E_{n,1} \cap E_{n,2} \cap E_{n,3})  | \leq c (\log n)^{ \frac{3d^2} {2} + \frac{d} {2}}
$$
or more simply
$$
 \E | ( F(  \hat{\P}^{( n, \alpha)} ) -  F(  \hat{\P}^{( n, \alpha)} \cup Z_{N(b_n) + 1} )) {\bf 1} ( E_{n,1} ) | \leq  (\log n)^{ \frac{3d^2} {2} + \frac{d} {2}}.
$$
Iterating this $j$ times gives
\begin{align*}
& \E | ( F(  \hat{\P}^{( n, \alpha)} ) -  F(  \hat{\P}^{( n, \alpha)} \cup \{ Z_{N(b_n) + 1},... Z_{N(b_n) + j} \}  )) {\bf 1} ( E_{n,1} ) | \\
& \leq cj (\log n)^{ \frac{3d^2} {2} + \frac{d} {2}}.
 \end{align*}
Let $j : = {\rm{Bi}}(n, \frac{b_n} {n}) - N(b_n)$ and  condition on the event that $j \leq \gamma \log n \sqrt{b_n}.$   On the complement of this event we  have $|{\rm{Bi}}(n, b_n/n) - b_n| \geq \gamma/2 \log n \sqrt{b_n}$ or $|N(b_n) - b_n| \geq \gamma/2 \log n \sqrt{b_n}$.  When $\gamma$ is large enough, standard tail bounds for the binomial and Poisson random variables show that the complement has probability $o(n^{-d})$ and thus
\begin{align*}
& \E | ( F(  \hat{\P}^{( n, \alpha)} ) -  F(  \hat{\P}^{( n, \alpha)} \cup \{ Z_{N(b_n) + 1},... Z_{N(b_n) + j} \}  )) {\bf 1} ( E_{n,1} ) | \\
& \leq c (\log n)^{ \frac{3d^2} {2} + \frac{d} {2} + 1}  \sqrt{b_n}
\end{align*}
i.e.,
$$
 \E | ( F(  \hat{\P}^{( n, \alpha)} ) -  F(  \hat{\X}^{( n, \alpha)}   )) {\bf 1} ( E_{n,1} ) | \leq c (\log n)^{ \frac{3d^2} {2} + \frac{d} {2} + 1}   \sqrt{b_n}.
$$
In view of the bound \eqref{boundfora} we deduce
$$
\E | ( F(  \hat{\P}^{( n, \alpha)} ) -  F(  \hat{\X}^{( n, \alpha)}   )) {\bf 1} ( E_{n,1} ) | = o(  u_{n,\alpha}^{d-1} ).
$$
We have the same bound on the event  $\{ N(b_n) \geq {\rm{Bi}}(n, \frac{b_n} {n} )\}$.
Thus \eqref{smalldifference} holds and
\be  \label{binomial-limit1-new}
\lim_{\la \to \infty} \frac{ \E F(  \hat{\X}^{( n, \alpha)} ) } { n }  =
 \int_{0}^{\infty}  \E \xi_k^{(\infty, \alpha)} ((\0,h), {\P}^{(\infty, \alpha)} ) \mathrm{d}\nu^{(\infty,\alpha)}(h).
\ee
The height bound \eqref{Height} also holds for binomial input.  Appealing to this bound, we find that \eqref{binomial-limit1-new} also holds when  $\hat{\X}^{( n, \alpha)}$ is replaced by  ${\X}^{( n, \alpha)}$.

\vskip.2cm

\noindent {\em (ii) The case $\alpha \in (- \infty, \frac{-2} {d-1}]$}.  As above, enumerate the points in $ {\cal P}^{(n, \alpha)}$ by $Z_1,...,Z_{N(n)}$. As above we consider the coupled point set ${\cal Y}_n$, with ${\cal Y}_n \eqd \X^{(n, \alpha)}$.
We do not restrict $ {\cal P}^{(n, \alpha)}$ to a subset of $W_n$.

For any $\X \subset W_n$, define
\be \label{basicsum-2}
F(\X) : = F^{(n, \alpha)}(\X) := \sum_{x \in  \X } \xi_k^{(n, \alpha)} (x, \X ).
\ee

It suffices  to show that
\be \label{smalldifference1}
\E |  F(  {\P}^{( n, \alpha)} ) -  F(  {\X}^{( n, \alpha)}   )  | = o(n).
\ee

We follow the proof for the case $\alpha \in (\frac{-2} {d-1}, \infty)$ and put $b_n = n$, ${\rm{Bi}}(n, \frac{b_n} {n}) = n$.
Put  $E_{n, 1} := \{ N(n) \leq n \}$.
Let $E_{n,2} $ be the event that the radii of stabilization of $\xi^{(n, \alpha)}$ at all points in  ${\P}^{( n, \alpha)}$ are bounded by $\gamma \log n$. If $\gamma$ is large enough then $\PP(E_{n,2} ^c) = o (n^{-d - 1})$.
Conditional on $Z_{N(n) + 1} = (v_0, h_0)$, on the event $E_{n,2}$, the number of extreme points in ${\P}^{( n, \alpha)}$ whose score may be modified by the insertion of $Z_{N(n) + 1}$ is bounded by the  number of points of ${\P}^{( n, \alpha)}$ with spatial coordinates
in $B_{d-1}(v_0, \gamma \log n)$.  There is an event $E_{n, 4}$ with $\PP(E_{n, 4}^c) = o(n^{-d})$ such that on $E_{n, 4}$ this number is bounded by $c (\log n)^{d-1}$.

This gives
$$
 \E | ( F(  {\P}^{( n, \alpha)} ) -  F(  {\P}^{( n, \alpha)} \cup Z_{N(n) + 1} )) {\bf 1} (E_{n,1} \cap E_{n,2} \cap E_{n,4})|
  \leq c (\log n)^{ \frac{3d^2} {2} + \frac{d} {2}}
$$
or more simply
$$
 \E | ( F(  {\P}^{( n, \alpha)} ) -  F(  {\P}^{( n, \alpha)} \cup Z_{N(n) + 1} )) {\bf 1} ( E_{n,1} ) | \leq c (\log n)^{ \frac{3d^2} {2} + \frac{d} {2}} .
$$
Iterating this $j$ times gives
\begin{align*}
&  \E | ( F(  {\P}^{( n, \alpha)} ) -  F(  {\P}^{( n, \alpha)} \cup \{ Z_{N(b_n) + 1},... Z_{N(n) + j} \}  )) {\bf 1} ( E_{n,1} ) | \\
& \leq cj (\log n)^{ \frac{3d^2} {2} + \frac{d} {2}} .
\end{align*}
Let $j : = n - N(n)$ and  condition on the event that $j \leq \gamma \log n \sqrt{n}.$    When $\gamma$ is large enough, standard tail bounds for the Poisson random variable show that the complement has probability $o(n^{-d})$ and thus
$$
\E | ( F(  {\P}^{( n, \alpha)} ) -  F(  {\X}^{( n, \alpha)}   )) {\bf 1} ( E_{n,1} ) | = o(n).
$$
We have the same bound on the event  $\{ N(n) \geq n \}$.
Thus \eqref{smalldifference1} holds as desired.
\qed

\vskip.2cm

\noindent{\bf 4.2. Proof of Theorem \ref{Thm2}.}
We only show the convergence for  Poisson initial input as the proof for the binomial input follows along similar lines.

The result is a consequence of Billingsley's continuous mapping theorem, see \cite[Theorem 5.5]{Bil}. The idea is to rewrite the set of rescaled extreme points as a continuous function of the initial input. To do so, we start by endowing the set ${\mathcal L}$ of locally finite sets of $\R^{d-1}\times \R^+$ with the distance
$\delta(\chi_1,\chi_2): =\sum_{n\ge 1}2^{-n} {\bf 1} (\chi_1\cap S({\bf 0},n,n)\neq \chi_2\cap S({\bf 0},n,n))$
 for all $\chi_1,\chi_2\in {\mathcal L}$, where $S({\bf 0},n,n):= B_{d-1}(0,n) \times [0,n].$
 Fix $\alpha \in \R$. For any $\chi\in {\mathcal L}$ and $\lambda \in [1, \infty)$, let $h_{\lambda}(\chi)$  be the subset $\rm{Ext}_\lambda(\chi)$ of $\chi$ comprised of points which survive the quasi parabolic thinning corresponding to the quasi parabolic festoon generated by $[ \Pi^{\downarrow}( \cdot) ]^{(\la, \alpha)}$ described in Section 3.1. This quasi parabolic festoon depends on the parameter $\alpha$ which implies that both $h_\lambda$ and $\rm{Ext}_\lambda$ depend on $\alpha$ but for sake of simplicity, we do not make this dependency visible in the notation.
Similarly, we define $h(\chi)$ to be the subset $\rm{Ext}(\chi)$ of $\chi$ generated by  $[ \Pi^{\downarrow}( \cdot) ]^{(\infty, \alpha)}$; see Section 1.3.

Let us show the continuity assumption which is required to apply the continuous mapping theorem: Let $(\chi_\lambda)$ be a parametrized family of elements of ${\mathcal L}$ which converges to $\chi$ when $\lambda\to\infty$.
%We assume that $\chi_\infty$ is in regular position, which happens with probability one for the considered Poisson point process. Then ecause of the topology provided by the distance $\delta$,
We have that for any fixed $n$, $\chi_\lambda$ and $\chi$ coincide on $S({\bf 0},n,n)$ for $\lambda$ large enough. Since the determination of the intersection of $\rm{Ext}(\chi)$ with $S({\bf 0},n,n)$ only depends on $\chi\cap S({\bf 0},m,m)$ for some $m\ge n$, we obtain that for any $n$,  $\rm{Ext}(\chi_\lambda)$ and $\rm{Ext}(\chi)$ coincide on $S({\bf 0},n,n)$  for $\lambda$ large enough.
%(Ext_\infty(\chi_\infty) may depend on points of \chi_\infty in S(m), m>n, but it then works as soon as \chi_\lambda and \chi_\infty coincide on S(m)).

Recall the definition of `regular position' of  point sets given in subsection \ref{expectation-1}.
We claim that once $\chi\cap S({\bf 0},n,n)$ is in regular position, $\rm{Ext}_\lambda(\chi_\lambda)$ and $\rm{Ext}(\chi_\lambda)$ coincide in $S({\bf 0},n,n)$ for $\lambda$ large enough.
Now our claim may be justified as follows:   As soon as the set of underlying points $\chi\cap S({\bf 0},n,n)$ is fixed, the rule depending on $\lambda$ to decide whether a point is extreme (i.e. to decide whether the point survives the quasi parabolic thinning), provides the same result for $\lambda$ large enough because the quasi-parabolic grains converge uniformly on any compact set to parabolic grains.
This implies that $h_\lambda(\chi_\lambda)$ converges to $h(\chi)$ when $\chi$ is in regular position.

Since a.s. the limiting Poisson point process $\P^{(\infty, \alpha)}$ is in regular position,  the continuity assumption of the continuous mapping theorem is satisfied  by taking  the set $E$ of  \cite[Theorem 5.5]{Bil} to be the set of locally finite point sets of $\R^{d-1}\times\R^+$ which are not in regular position. Applying that theorem, we get that  $\rm{Ext}_\lambda(\P^{(\lambda, \alpha)} )$ converges in distribution as $\la \to \infty$
to $\rm{Ext}(\P^{(\infty, \alpha) })$, as required. \qed

\vskip.2cm

\noindent{\bf 4.3. Proof of Theorem \ref{variance}. }
Before proving Theorem \ref{variance} we need one lemma. We extend the definition of the correlation function at \eqref{SO2aa} and define for all $\la \in [1, \infty]$, $w_1, w_2 \in W_\la$ and $\alpha \in (- \infty, \infty)$
\begin{align*}
c^{\xi_k^{(\la, \alpha)} }(w_1, w_2) & := \E \xi_k^{(\la, \alpha)}(w_1, {\P}^{(\la, \alpha)} \cup \{w_2 \} ) \xi_k^{(\la, \alpha)}(w_2, {\P}^{(\la, \alpha)} \cup \{w_1 \} )\\
& \ \ \ \ \ \ \ \ \ -  \E \xi_k^{(\la, \alpha)}(w_1, {\P}^{(\la, \alpha)} ) \E \xi_k^{(\la, \alpha)}(w_2, {\P}^{(\la, \alpha)} ).
\end{align*}
Proposition \ref{expectation} implies for all $\alpha \in [\frac{-2} {d-1}, \infty)$ and all $w_1, w_2 \in W_\la$ that
\be \label{2ptconvergence}
\lim_{\la \to \infty} c^{\xi_k^{(\la, \alpha)} }(w_1, w_2) = c^{\xi_k^{(\infty, \alpha)} }(w_1, w_2).
\ee
On the other hand, for all $\alpha \in (-\infty, \frac{-2} {d-1} )$ and all $(v_1, h_1), (v_2, h_2) \in W_\la$,
Proposition \ref{expectation-1} gives
\be \label{2ptconvergence-1}
\lim_{\la \to \infty} c^{\xi_k^{(\la, \alpha)} }((v_1, h_1), (v_2, h_2) ) = c^{\xi_k^{(\infty, \alpha)} }((v_1, 0), (v_2, 0)).
\ee

 The two point correlation function $c^{\xi_k^{(\la, \alpha)} }(w_1, w_2)$ decays exponentially fast with the distance between the spatial coordinates of $w_1$ and $w_2$.

%\Comment{part (ii) of Lemma \ref{covariance-decay} is new}
\begin{lemm} \label{covariance-decay}  (i) For all $\alpha \in (\frac{-2} {d-1}, \infty)$ there is a constant $c \in (0, \infty)$ such that for all $w_1:= (v_1, h_1)$ and $w_2:= (v_2, h_2)$ in $W_\la, \la \in [1, \infty],$ we have
\be \label{covdecay}
c^{\xi_k^{(\la, \alpha)} }(w_1, w_2) \leq c (h_1 h_2)^{c} \exp( - \frac{ h_1^{\frac{d + 1} {2} } + h_2^{ \frac{d + 1} {2} } } {c} )  \exp( - \frac{ |v_1 - v_2|^{d + 1} } {c} ).
\ee
(ii) When  $\alpha  \in (- \infty,  \frac{-2} {d-1}]$ there is a constant $c \in (0, \infty)$ such that for all $w_1:= (v_1, h_1)$ and $w_2:= (v_2, h_2)$ in $W_\la, \la \in [1, \infty],$ we have
\be \label{covdecay-1}
c^{\xi_k^{(\la, \alpha)} }(w_1, w_2) \leq c  \exp( - \frac{ h_1^{  \frac{d - 1} {2} } + h_2^{\frac{ d - 1} {2}} } {c} )  \exp( - \frac{ |v_1 - v_2|^{d-1} } {c} ).
\ee
\end{lemm}

\noindent{\em Proof.} (i) Abbreviate $\xi_k^{(\la, \alpha)}$ by $\xi$. Let $p > 2$.
  Put $R := \max (R_{w_1}, R_{w_2})$, where $$R_{w_1}:= R^{\xi^{(\la, \alpha)}}(w_1, {\P}^{(\la, \alpha)}),  \ \ R_{w_2}:= R^{\xi^{(\la, \alpha)}}(w_2, {\P}^{(\la, \alpha)})$$ are the radii of stabilization  at $w_1$ and $w_2$, respectively. Furthermore, put $r:= \frac{|v_1 - v_2 |} {3}$ and define the event $E := \{ R \in [0, r] \}$.  H\"older's inequality gives
\begin{align} \label{cov1}
& | \E \xi(w_1, {\P}^{(\la, \alpha)}  \cup \{w_2\} ) \xi(w_2, {\P}^{(\la, \alpha)} \cup \{w_1\} ) \nonumber \\
& \ \ \ \ \ \ \ \ \  - \E \xi(w_1, {\P}^{(\la, \alpha)} \cup \{w_2\} ) \xi(w_2, {\P}^{(\la, \alpha)} \cup \{w_1\} )\1(E) | \nonumber \\
& \leq c \left(\sup_{  w_1, w_2 \in \R^{d-1} \times \R^+   } \E | \xi(w_1, {\P}^{(\la, \alpha)} \cup \{w_2\} )|^p\right)^{ \frac{2} {p} } \, \PP(E^c)^{ \frac{p-2} {p}}.
\end{align}
Notice that
\begin{align*}
& \E \xi(w_1, {\P}^{(\la, \alpha)} \cup \{w_2\} ) \xi(w_2, {\P}^{(\la, \alpha)} \cup \{w_1\} ) \1(E)  \\
& = \E \xi(w_1, ({\P}^{(\la, \alpha)}  \cup \{w_2\}) \cap C_{d-1}(v_1,r)  ) \\
& \ \ \ \  \ \ \ \ \ \ \times \xi(w_2, ({\P}^{(\la, \alpha)} \cup \{w_1\}) \cap C_{d-1}(v_2,r)  ) \1(E) \\
& = \E \xi(w_1, ({\P}^{(\la, \alpha)} \cup \{w_2\}) \cap C_{d-1}(v_1,r)  ) \\
& \ \ \ \ \ \ \ \ \ \ \times \xi(w_2, ({\P}^{(\la, \alpha)} \cup \{w_1\}) \cap C_{d-1}(v_2,r)  ) (1 -  \1(E^c)).
\end{align*}

A second application of  H\"older's inequality gives
\begin{align} \label{cov2}
& | \E \xi(w_1, {\P}^{(\la, \alpha)} \cup \{w_2\} ) \xi(w_2,{\P}^{(\la, \alpha)} \cup \{w_1\} ) \1(E)   \nonumber \\
& \ \ \ \ \ \  - \E \xi(w_1, ( {\P}^{(\la, \alpha)} \cup \{w_2\})  \cap C_{d-1}(v_1,r) ) \xi(w_2, ({\P}^{(\la, \alpha)} \cup \{w_1\})  \cap C_{d-1}(v_2,r) )  | \nonumber \\
& \leq c \left(\sup_{  w_1, w_2 \in \R^{d-1} \times \R^+  } \E | \xi(w_1, {\P}^{(\la, \alpha)} \cup \{w_2\} )|^p\right)^{ \frac{2} {p} } \,
\PP(E^c)^{ \frac{p-2} {p} }.
\end{align}
\vskip.2cm
Combining \eqref{cov1} and \eqref{cov2} and using that  $\xi(w_1, ({\P}^{(\la, \alpha)}  \cup \{w_2\})  \cap C_{d-1}(v_1,r) )$ and
$\xi(w_2, ({\P}^{(\la, \alpha)}  \cup \{w_1\})  \cap C_{d-1}(v_2,r) )$ are independent we have
\begin{align*} %\label{cov3}
& | \E \xi(w_1, {\P}^{(\la, \alpha)} \cup \{w_2\} ) \xi(w_2, {\P}^{(\la, \alpha)} \cup \{w_1\} )\nonumber \\
& \ \ \ - \E \xi(w_1, ({\P}^{(\la, \alpha)}  \cup \{w_2\})  \cap C_{d-1}(v_1,r) ) \E \xi(w_2, ({\P}^{(\la, \alpha)} \cup \{w_1\})  \cap C_{d-1}(v_2,r) )  | \nonumber  \\
& \leq c \left(\sup_{  w_1, w_2 \in \R^{d-1} \times \R^+    } \E | \xi(w_1, {\P}^{(\la, \alpha)} \cup \{w_2\} )|^p\right)^{ \frac{2} {p} } \, \PP(E^c)^{ \frac{p-2} {p} }.
\end{align*}
\vskip.1cm
Likewise we may show
\begin{align*} % \label{cov4}
& | \E \xi(w_1, {\P}^{(\la, \alpha)}  ) \E \xi(w_2, {\P}^{(\la, \alpha)}  ) \nonumber \\
& \ \ \ \ \ \  - \E \xi(w_1, {\P}^{(\la, \alpha)}  \cap C_{d-1}(v_1,r) ) \E \xi(w_2, {\P}^{(\la, \alpha)} \cap C_{d-1}(v_2,r) )  |  \nonumber \\
& \leq c \left(\sup_{  w_1, w_2 \in \R^{d-1} \times \R^+    } \E | \xi(w_1, {\P}^{(\la, \alpha)} \cup \{w_2\} )|^p\right)^{ \frac{2} {p} } \, \PP(E^c)^{ \frac{p-2} {p} }.
\end{align*}
Combining the last two displays with %\eqref{cov3}, \eqref{cov4},
Lemma \ref{mombounds}, and using that $\PP(E^c)$ decays exponentially with $|v_1 - v_2|^{d + 1}$, we obtain \eqref{covdecay}.

\vskip.2cm
\noindent (ii) Follow the  proof for part (i), using instead the moment bounds of Lemma \ref{mombounds}(ii).
\qed

\vskip.2cm

\noindent{\em Proof of Theorem \ref{variance}. } By the Slivnyak - Mecke formula we have
\begin{align*}
\Var f_k(K_{\la,\alpha}) & = \int_{B(\0, 1 + \la^{\alpha})}  \E \xi_k^2(x, \tilde{\P}_{\la, \alpha}) \mathrm{d} \tilde{\P}_{\la, \alpha}(x) \\
&  \ \ \ \ + \int_{B(\0, 1 + \la^{\alpha})}  \int_{B(\0, 1 + \la^{\alpha})} [ \E \xi_k(x, \tilde{\P}_{\la, \alpha} \cup \{y \})  \xi_k(y, \tilde{\P}_{\la, \alpha} \cup \{x \}) \\
& \ \ \ \ \ \ \ - \E \xi_k(x, \tilde{\P}_{\la, \alpha})\E \xi_k(y, \tilde{\P}_{\la, \alpha}) ] \mathrm{d} \tilde{\P}_{\la, \alpha} (x)\mathrm{d} \tilde{\P}_{\la, \alpha} (y) \\
& := J_1(\la) + J_2(\la).
\end{align*}

 As in the proof of  Theorem \ref{expectasy} we deduce that
\be \label{J1int}
\lim_{\la \to \infty} \frac{ J_1(\la) } {d\kappa_d u_{\la,\alpha}^{d-1}  }  =  \int_{0}^{\infty}  \E (\xi_k^{(\infty, \alpha)} ((\0,h), {\P}^{(\infty, \alpha)} ))^2  \mathrm{d}\nu^{(\infty,\alpha)}(h).
\ee
Now we establish the convergence of $J_2(\la)/d \kappa_d u_{\la,\alpha}^{d-1}$ as $\la \to \infty$.
Given ${\bf u} \in \S^{d-1}$, ${\bf u} \neq {\bf u}_0$,
we define $T^{(\la, \alpha)} $ exactly as in \eqref{map2}, but with ${\bf u}_0:= (0,0,...,1) \in \R^d$ replaced by ${\bf u}$. We write $T_{\bf u}^{(\la, \alpha)}$ to denote the dependency on ${\bf u}$.

Denoting by $(\0, h_0)$ and $(v_1, h_1)$ the images under $T_{x/|x|}^{(\la, \alpha)}$ of $x$ and $y$, respectively, we notice that the `covariance' term in the integrand of $J_2(\la)$ transforms to the covariance $c^{\xi_k^{(\la, \alpha)} }((\0, h_0), (v_1, h_1) ).$  We apply the following change of variables in the quadruple integral:
$$
{\bf u} =\frac{x}{|x|},\quad h_0 = u_{\la,\alpha}^2(1-\frac{|x|} {1+ \la^{\alpha}}), \quad (v_1,h_1)=T_u^{(\la,\alpha)}(y).
$$
The double integral over $B(\0, 1 + \la^{\alpha}) \times B(\0, 1 + \la^{\alpha})$ transforms into a quadruple integral over $\S^{d-1} \times [0, h^{(\la,\alpha)}] \times u_{\la,\alpha} B_{d-1}(\pi) \times [0, h^{(\la,\alpha)} ].$   The intensity measure $\mathrm{d} \tilde{\P}_{\la, \alpha}(y)$ transforms to $\mathrm{d} \mu^{(\la, \alpha)}$ while
%, which converges to $\mathrm{d}\mu^{(\infty,\alpha)}(v,h)$ as $\la \to \infty$.
the intensity measure $\mathrm{d} \tilde{\P}_{\la, \alpha}(x)$ transforms to the
product measure $ u_{\la,\alpha}^{d-1} \mathrm{d}\sigma_{d-1}({\bf u}) \varphi^{(\la,\alpha)}(h)\mathrm{d}h.$

 In other words the term $J_2(\la)/ d\kappa_d u_{\la,\alpha}^{d -1}$ transforms to
 \begin{alignat}{5}\label{eq:rewritingJ2useful}
 & \int_0^{h^{(\la,\alpha)}}
\int_{ u_{\la,\alpha} B_{d-1}(\pi) \times [0, h^{(\la,\alpha)}] }
c^{\xi_k^{(\la, \alpha)} }((\0, h_0), (v_1, h_1) )\nonumber\\&\hspace*{3.0cm}
\times \frac{\sin^{d-2}(u_{\la,\alpha}^{-1}|v_1|)}{|u_{\la,\alpha}^{-1}v_1|^{d-2}}\mathrm{d}v_1\varphi^{(\la,\alpha)}(h_1)\mathrm{d}h_1
\varphi^{(\lambda,\alpha)}(h_0)\mathrm{d}h_0.
 \end{alignat}

\noindent (i) {\em The case $\alpha \in [\frac{-2} {d-1}, \infty)$} .
As in the proof of Theorem \ref{expectasy}, truncating the integrals in the height coordinates $h_0$ or $h_1$  at level $\log \la$ induces negligible error.  In other words,
\begin{align*}
\frac{ J_2(\la) }  { d\kappa_du_{\la,\alpha}^{d -1}  }  & =
\int_0^{\log \la}
\int_{ u_{\la,\alpha} B_{d-1}(\pi) \times [0, \log\la] }
c^{\xi_k^{(\la, \alpha)} }((\0, h_0), (v_1, h_1) )
\\ & \ \ \ \ \ \ \ \times \frac{\sin^{d-2}(u_{\la,\alpha}^{-1}|v_1|)}{|u_{\la,\alpha}^{-1}v_1|^{d-2}}\mathrm{d}v_1\varphi^{(\la,\alpha)}(h_1)\mathrm{d}h_1
\varphi^{(\lambda,\alpha)}(h_0)\mathrm{d}h_0+o(1)
\end{align*}
where we have also used Lemma \ref{lem:proppn} (i).

By \eqref{eq:convphi} and Proposition \ref{expectation}, the integrand above converges a.e. to
$$c^{\xi_k^{(\la, \alpha)} }((\0, h_0), (v_1, h_1) )\varphi^{(\infty,\alpha)}(h_1)\varphi^{(\infty,\alpha)}(h_0)
%\mathrm{d}\mu^{(\infty,\alpha)}(v_1,h_1)\mathrm{d}\mu_h^{(\infty,\alpha)}(h_0)
.$$
%$$c^{\xi_k^{(\la, \alpha)} }((\0, h_0), (v_1, h_1) ) d {\mu}^{(\la, \alpha)}(v_1,h_1)  \varphi^{(\lambda,\alpha)}(h_0)\mathrm{d}h_0$$
By Lemma \ref{covariance-decay} and \eqref{eq:upperboundphi}, the above is dominated uniformly in $\la$ by a function decaying exponentially fast in $h_0, h_1$, and in $|v_1|^{d + 1}$ on the region $[0, \log \la] \times u_{\la,\alpha} B_{d-1}(\pi) \times [0, \log \la].$   Applying the dominated convergence theorem,
we obtain
 \begin{align*}
 & \lim_{\la \to \infty} \frac{ J_2(\la) }  { d\kappa_du_{\la,\alpha}^{d -1} } \\
 &  =   \int_{0}^{\infty}  \int_{\R^{d-1}} \int_{0}^{\infty}
   c^{\xi_k^{(\infty, \alpha )}}((\0,h_0),(v_1,h_1))\varphi^{(\infty,\alpha)}(h_1)\varphi^{(\infty,\alpha)}(h_0) \mathrm{d}h_1 \mathrm{d}v_1\mathrm{d} h_0\\
 &  =   \int_{0}^{\infty}  \int_{\R^{d-1}} \int_{0}^{\infty}
   c^{\xi_k^{(\infty, \alpha )}}((\0,h_0),(v_1,h_1)) \mathrm{d} \mu^{(\infty,\alpha)} (v_1,h_1) \mathrm{d} \nu^{(\infty,\alpha)}(h_0).
  \end{align*}
Thus
$$
\lim_{\la \to \infty} \frac{ J_1(\la) + J_2(\la) }  { d\kappa_du_{\la,\alpha}^{d -1} } = \sigma^2(\xi_k^{(\infty, \alpha)})
$$
which is the desired convergence as at \eqref{lim1}.
Putting $\la = \infty$ in  Lemma \ref{covariance-decay}(i), we conclude that
$c^{\xi_k^{(\infty, \alpha )}}$ is integrable and thus $\sigma^2(\xi_k^{(\infty, \alpha)}) < \infty$.  The proof that $\sigma^2(\xi_k^{(\infty, \alpha)})$ is strictly positive is postponed to subsection 4.5.

\vskip.2cm

\noindent (ii) {\em The case $\alpha \in (-\infty, \frac{-2} {d-1})$.}
Let us consider $\lim_{\la \to \infty} J_2(\la)/\la$. In view of the dominated convergence theorem, we have
\begin{align*}
&\int_{\R^{d-1}}  c^{\xi_k^{(\infty, \alpha)}}((\0,0),(v_1,0)) \mathrm{d}v_1\nonumber \\ &=\int_{u_{\la,\alpha} B_{d-1}(\pi)} c^{\xi_k^{(\infty, \alpha)}}((\0,0),(v_1,0))\frac{\sin^{d-2}(u_{\la,\alpha}^{-1}|v_1|)}{|u_{\la,\alpha}^{-1}v_1|^{d-2}} \mathrm{d}v_1+o(1).
\end{align*}
Integrating the right-hand side  with respect to $\varphi^{(\la,\alpha)}(h_1)\varphi^{(\la,\alpha)}(h_0)\mathrm{d}h_1\mathrm{d}h_0$ and using \eqref{eq:totalmassmulambdaalpha}, we get
\begin{align}
  \label{eq:rewritinglimitvar2}
& \int_{\R^{d-1}}  c^{\xi_k^{(\infty, \alpha)}}((\0,0),(v_1,0)) \mathrm{d}v_1\nonumber \\
& =\int_0^{h^{(\la,\alpha)}}\int_0^{h^{(\la,\alpha)}}\int_{u_{\la,\alpha} B_{d-1}(\pi)} c^{\xi_k^{(\infty, \alpha)}}((\0,0),(v_1,0))\nonumber\\
& \ \ \ \ \ \ \times \frac{\sin^{d-2}(u_{\la,\alpha}^{-1}|v_1|)}{|u_{\la,\alpha}^{-1}v_1|^{d-2}}\varphi^{(\la,\alpha)}(h_0)\varphi^{(\la,\alpha)}(h_1)  \mathrm{d}v_1\mathrm{d}h_1\mathrm{d}h_0 +o(1).
\end{align}
Comparing \eqref{eq:rewritinglimitvar2} with \eqref{eq:rewritingJ2useful}, we deduce that it is enough to prove that the integral with respect to $v_1$ of
\begin{align}
  \label{eq:integrandfunction}
&\int_0^{h^{(\la,\alpha)}}\int_0^{h^{(\la,\alpha)}}| c^{\xi_k^{(\la, \alpha)} }((\0, h_0), (v_1, h_1) ) -  c^{\xi_k^{(\infty, \alpha )}}((\0, 0),(v_1,0))|\nonumber\\&\hspace*{3cm}\times \frac{\sin^{d-2}(u_{\la,\alpha}^{-1}|v_1|)}{|u_{\la,\alpha}^{-1}v_1|^{d-2}}\varphi^{(\la,\alpha)}(h_0)\varphi^{(\la,\alpha)}(h_1)\mathrm{d}h_0\mathrm{d}h_1
\end{align}
converges to zero. To do so, we apply the dominated convergence theorem to this simple integral in $v_1$.
First, recalling \eqref{LconvPo-1} and \eqref{LconvPo-3} and letting $\la \to \infty$,  we obtain
$$ %\label{LconvPo-4}
\sup_{h_0\in [0,h^{(\lambda,\alpha)}]} \sup_{(v_1,h_1) \in W_\la}
| c^{\xi_k^{(\la, \alpha)} }((\0, h_0), (v_1, h_1) ) -  c^{\xi_k^{(\infty, \alpha )}}((\0, 0),(v_1,0))| \to 0.
$$
This implies, thanks to \eqref{eq:totalmassmulambdaalpha}, that the double integral at \eqref{eq:integrandfunction} converges to zero for all $v_1$.
Secondly, thanks to Lemma \ref{covariance-decay} and again \eqref{eq:totalmassmulambdaalpha}, we notice that this same quantity is dominated by an exponentially decaying integrable function in $v_1$. Consequently, we get
$$
\lim_{\la \to \infty} \frac{ J_{2}(\la)  }  {\la} =   \int_{ \R^{d-1}  } c^{\xi_k^{(\infty, \alpha)} }((\0, 0), (v_1, 0) ) \mathrm{d}v_1.
$$
Together with \eqref{J1int}, we obtain
\begin{align*}
& \lim_{\la \to \infty} \frac{  J_{1}(\la) +   J_{2}(\la)  }  { \la  } \\
&  =   \E ( \xi_k^{(\infty, \alpha)} ((\0,0), {\P}^{(\infty, \alpha)} ))^2
+   \int_{ \R^{d-1}  } c^{\xi_k^{(\infty, \alpha)} }((\0, 0), (v_1, 0) ) \mathrm{d}v_1 \\
& = \sigma^2(\xi_k^{(\infty, \alpha)}).
 \end{align*}
This is the desired convergence  \eqref{S03}. Note that $\sigma^2(\xi_k^{(\infty, \alpha)})$ is finite by Lemma \ref{covariance-decay}(ii). This shows
\eqref{lim1} when $\alpha \in (-\infty, \frac{-2} {d-1})$.  The validity of \eqref{lim1a} holds since  $c^{\xi_0^{(\infty, \alpha)} }((\0, 0), (v_1, 0) )$ vanishes and $\xi_0^{(\infty, \alpha)} ((\0, 0), \P^{(\infty, \alpha)}) = 1$.
This concludes the proof of \eqref{lim1} when $\alpha \in (- \infty, \frac{-2} {d-1})$, save for showing that $\sigma^2(\xi_k^{(\infty, \alpha)})$ is positive.  This will be shown in subsection 4.5.  \qed

\vskip.2cm

\noindent{\bf 4.4. Proof of Theorem \ref{clt}.}  We follow the proof of Theorem 7.1 of \cite{CSY}, which is based on dependency graph arguments.  Our proof is simpler since we do not show a central limit theorem for random measures, but simply for their total mass.  In other words, we  put the function $g$ in Theorem 7.1 of \cite{CSY} to be identically one. We also put $\delta = 0$ in that theorem.  When the value of the parameter $\alpha$ exceeds $\frac{2} {d +1}$, the proof is unchanged since the scaling transform is identical.  For the other values of the parameter $\alpha$, it suffices to  replace $\beta$ by $\beta(\alpha)$ and to make the identification $\tau = (d -1) \beta(\alpha)$.  \qed

\vskip.2cm

\noindent{\bf 4.5. Positivity of the limiting variance  $\sigma^2(  \xi_k^{(\infty, \alpha)})$.}

\noindent {\em (i) The case $\alpha \in [\frac{-2} {d-1}, \infty)$.}
Define
$$
H_{k, \lambda}^{(\infty, \alpha)} := \sum_{x \in \P^{(\infty, \alpha)} \cap W_\la } \xi_k^{(\infty, \alpha)} (x, \P^{(\infty, \alpha) }).
$$
The proof of variance asymptotics \eqref{lim1} for  the case $\alpha \in [\frac{-2} {d-1}, \infty)$
is easily adapted to show that $\sigma^2(  \xi_k^{(\infty, \alpha)})$ coincides with
$$
\lim_{\la \to \infty} \frac{ \Var H_{k, \lambda}^{(\infty, \alpha)}} {\la^{\beta(d -1)} }.
$$
Thus it suffices to show that the above limit is positive.

Our approach is based on two observations. First, let $Q$ be a fixed cube in $\R^{d-1}$. Fx $\rho \in (0, \infty)$ with $\rho$ smaller than the diameter of $Q$.  Consider the following deterministic union of paraboloids:  At each point $x_0$ of the grid $\rho \Z^{d-1} \cap Q$, we consider the `up-paraboloid' $y \geq   ||x - x_0||^2/2$ and we take the union over all such paraboloids having apices in $\rho \Z^{d-1} \cap Q$.  The boundary of this union is called a parabolic hull. Denote the number of $k$-faces, $k \in \{0,1,...,d-1 \}$, of the parabolic hull by $F_k(Q, \rho)$.  All vertices of the parabolic hull are extreme points. Scaling gives $F_0(Q, \frac{\rho} {2}) = 2^{d-1} F_0(Q, \rho)$.  For more general $k$, this scaling relation
is not exact, as boundary effects  may play a role.  These effects are negligible for small $\rho$ and thus as $\rho \to 0$ we have
$$
F_k(Q, \frac{\rho} {2}) \sim 2^{d-1} F_k(Q, \rho),  \ k \in \{0,1,...,d-1 \}.
$$
Given $\rho$, there is a small positive $\delta := \delta(\rho),  \delta < < \rho,$ such that if the apices located at the grid points in $\rho \Z^{d-1} \cap Q$ are perturbed by at most $\delta$,  then the  number of $k$-faces, $k \in \{0,1,...,d-1 \}$, of the resulting $\delta$-modified parabolic hull, does not depend on the locations of the perturbed points.  We denote this number by $F_k(Q, \rho, \delta)$ and note that it coincides with
$F_k(Q, \rho)$. Similarly, we define $F_k(Q, \frac{\rho} {2}, \delta)$ and note that  as $\rho \to 0$ we have
\be \label{scaling}
F_k(Q, \frac{\rho} {2}, \delta) \sim 2^{d-1} F_k(Q, \rho, \delta),  \ k \in \{0,1,...,d-1 \}.
\ee

A second observation goes as follows. Given any closed set $B \subset \R^{d-1}$ and $\epsilon > 0$, define $B^{(\epsilon)}:= \{ x \in B: \ d(x, \partial B) > \epsilon \}.$
Let $\rho$ and $\delta := \delta(\rho)$ be as above. If the point process $\P^{(\infty, \alpha)}$ puts exactly one point in the {\em $d$-dimensional ball} of radius $\delta$ centered at each grid point in $\rho \Z^{d-1} \cap Q$, and if $Q \times [0, \delta]$ contains no other point from $\P^{(\infty, \alpha)}$, then we say that $\P^{(\infty, \alpha)}$ is an `approximate $\rho$-grid point process on $Q \times [0, \delta]$'.
 Given such a point process, observe that $F_k(Q^{(2\rho)}, \rho, \delta)$ is independent of the configuration of points in $\P^{(\infty, \alpha)} \cap (Q^c\times [0,\infty))$. Indeed,  when $\delta$ is small enough, the hull process defined by points in an `approximate $\rho$-grid point process on $Q^{(2\rho)} \times [0, \delta]$' will consist of paraboloids with `height' at most $\delta$ and thus
these paraboloids will never be entirely covered by paraboloids generated by points
of $\P^{(\infty, \alpha)}$
in $Q^c \times [0, \infty)$ (nor will the $k$-dimensional faces containing the apices of these paraboloids be completely covered).
We call this the `independence property' of the $F_k$ functional on the approximate $\rho$-grid point process.

Fix $\rho \equiv 1$ and take $\delta := \delta(1)$ as above.  Abusing notation we put $\tilde{W}_\la := [ \frac{ - \la^{\beta}} {2},  \frac{ \la^{\beta}} {2}]^{d-1}.$
Let $M$ be a large positive number. Partition $\tilde{W}_\la$  into a collection of  $L:= ( [ \frac{\la^{\beta}} {M} ])^{d-1}$ sub-cubes $Q_1,...,Q_L$.
Consider the sub-cubes $Q_i$ in this collection satisfying the following three properties:
\vskip.2cm
\noindent (a) $\P^{(\infty, \alpha)} \cap B_d(z, \delta)$ consists of a singleton for all $z \in \Z^{d-1} \cap ((Q_i \setminus Q_i^{(2)}) \times [0, \delta])$
\vskip.2cm
\noindent  (b) one of the following two events holds:
\vskip.1cm
(1) $\P^{(\infty, \alpha)}  \cap B_d(z, \delta)$ consists of a singleton for all $z \in \Z^{d-1} \cap Q_i^{(2)}$, or
\vskip.1cm
(2) $\P^{(\infty, \alpha)} \cap B_d(z, \delta)$ consists of a singleton for all $z \in \frac{1} {2} \Z^{d-1} \cap Q_i^{(2)}$
\vskip.2cm
\noindent  (c) $\P^{(\infty, \alpha)}$ puts no other point in $Q_i^{(2)} \times [0, \delta]$.
\vskip.2cm
Conditions (b1) and (c) imply that $\P^{(\infty, \alpha)}$ is an approximate $1$-grid point process on $Q_i^{(2)} \times [0, \delta]$, whereas conditions  (b2) and (c) imply that $\P^{(\infty, \alpha)}$ is an approximate $\frac{1} {2}$-grid point process on $Q_i^{(2)} \times [0, \delta]$.

Re-labeling if necessary, we let $I:= \{1,...,K \}$ be the indices of cubes having properties (a)-(c).  Elementary properties of the Poisson point process show that there is positive probability that a given cube satisfies conditions (a)-(c).  Thus $\E K = \Theta(\la^{\beta(d-1)}).$   Put ${\cal Q} := \cup_{i = 1}^K Q_i, {\cal Q}^c := \tilde{W}_\la \setminus {\cal Q}.$

Let ${\cal F}_\la$ be the $\sigma$-algebra determined by the random set $I$ and the positions of points of $\P^{(\infty, \alpha)}$ not belonging to $\cup_{i \in I} Q_i^{(2)} \times [0, \delta]$.
Given ${\cal F}_\la$, the functional $H_{k, \lambda}^{(\infty, \alpha)}$ admits variability only inside $\cup_{i \in I} Q_i^{(2)} \times [0, \delta]$.  Also, for each $i \in I$ we have
\be \label{lowerbound}
\Var [  \sum_{x \in \P^{(\infty, \alpha)}  \cap (Q_i^{(2)} \times [0, \delta])  }   \xi_k^{(\infty, \alpha)} (x, \P^{(\infty, \alpha) })  | {\cal F}_\la  ] \geq c_0
\ee
because conditional on ${\cal F}_\la$ either conditions $(b1)$ or $(b2)$ occur and because the scaling relation \eqref{scaling} gives two different values depending on whether $(b1)$ or $(b2)$ occurs.

For any random variable $X$ and $\sigma$-algebra $\cal F$, the conditional variance formula says that
$\Var X = \Var [ \E[X | {\cal F} ]] + \E [\Var[X | {\cal F} ]].$
Thus,
\begin{align*}
\Var [ H_{k, \lambda}^{(\infty, \alpha)}]  &  \geq \E  \Var [ H_{k, \lambda}^{(\infty, \alpha)}  | {\cal F}_\la  ] \\
& = \E \Var \left[ \sum_{i \in I} \sum_{x \in \P^{(\infty, \alpha)}  \cap (Q_i^{(2)} \times [0, \delta])  }   \xi_k^{(\infty, \alpha)} (x, \P^{(\infty, \alpha) })  | {\cal F}_\la  \right] \\
& = \E \sum_{i \in I} \Var \left[  \sum_{x \in \P^{(\infty, \alpha)}  \cap (Q_i^{(2)} \times [0, \delta])  }   \xi_k^{(\infty, \alpha)} (x, \P^{(\infty, \alpha) })  | {\cal F}_\la  \right] \\
& \geq c_0 \E K \\
& = \Theta(\la^{\beta(d-1)} ),
\end{align*}
where the second equality follows from the independence property of the
$F_k$ functional on $\P^{(\infty, \alpha)}$ on each $Q_i, i \in I,$ and where the last inequality follows from \eqref{lowerbound} applied to each $i \in I$.
Thus $\lim_{\la \to \infty}  \Var H_{k, \lambda}^{(\infty, \alpha)}/\la^{\beta(d -1)} $ is strictly positive, which was to be shown.

\vskip.2cm

\noindent {\em (ii) The case $\alpha \in (- \infty, \frac{-2} {d-1})$}.  We follow the same approach as for case (i).  Since $\P^{(\infty, \alpha)}$ is hosted by $\R^{d-1}$, we may restrict attention  to $\R^{d-1}$, making the arguments simpler.
We replace $B_d(z, \delta)$ by $B_{d-1}(z, \delta)$ and consider sub-cubes $Q_i$ having the properties:
\vskip.2cm
\noindent (a') $\P^{(\infty, \alpha)} \cap B_{d-1}(z, \delta)$ consists of a singleton for all $z \in \Z^{d-1} \cap (Q_i \setminus Q_i^{(2)})$
\vskip.2cm
\noindent  (b') one of the following two events holds:
\vskip.2cm
(1) $\P^{(\infty, \alpha)}  \cap B_{d-1}(z, \delta)$ consists of a singleton for all $z \in \Z^{d-1} \cap Q_i^{(2)}$, or
\vskip.2cm
(2) $\P^{(\infty, \alpha)} \cap B_{d-1}(z, \delta)$ consists of a singleton for all $z \in \frac{1} {2} \Z^{d-1} \cap Q_i^{(2)}$,
\vskip.2cm
\noindent  (c') $\P^{(\infty, \alpha)}$ puts no other points in $Q_i^{(2)}$.
\vskip.2cm
We let $I:= \{1,...,K \}$ be the indices of cubes having properties (a')-(c'). Let ${\cal F}_\la$ be the $\sigma$-algebra determined by the random set $I$, the positions of all points of $\P^{(\infty, \alpha)}$ not belonging to $\cup_{i \in I} Q_i^{(2)}$.
Now it suffices to follow the above arguments mutatis mutandis.
\qed

\vskip.2cm
\noindent{\em Acknowledgments}. A significant part of this research was completed at the Universit\'e de Rouen Normandie.   J. Yukich is grateful to the
department of mathematics for its kind hospitality and support.

\end{document}